\documentclass[a4paper]{amsart}
\RequirePackage{amsmath}
\RequirePackage{bm}
\RequirePackage{amssymb}
\RequirePackage{upref}
\RequirePackage{amsthm}
\RequirePackage{enumerate}
\RequirePackage{pb-diagram}
\RequirePackage{amsfonts}
\RequirePackage[mathscr]{eucal}
\RequirePackage{verbatim}
\RequirePackage{xr}
\RequirePackage{graphicx}
\usepackage{calc}
\usepackage{xspace}
\RequirePackage{color}
\RequirePackage{ifthen}

\newcommand{\cf}{cf.\@\xspace}
\newcommand{\resp}{resp.\@\xspace}

\newcommand{\al}{\alpha}
\newcommand{\bet}{\beta}
\newcommand{\ga}{\gamma}
\newcommand{\de}{\delta }
\newcommand{\e}{\epsilon}

\newcommand{\f}{\varphi}
\newcommand{\h}{\eta}

\newcommand{\ka}{\kappa}
\newcommand{\lam}{\lambda}

\newcommand{\n}{\nu}

\newcommand{\vt}{\vartheta}

\newcommand{\s}{\sigma}
\newcommand{\x}{\xi}

\newcommand{\C}{\varGamma}

\newcommand{\F}{\varPhi}
\newcommand{\Lam}{\varLambda}

\newcommand{\di}[1]{#1\nobreakdash-\hspace{0pt}dimensional}%\di n

\newcommand{\fmo}[1]{F_{|_{#1}}}%\fmo M

\newcommand{\fv}[2]{#1\hspace{0pt}_{|_{#2}}}

\newcommand{\const}{\tup{const}}

\newcommand{\msp[1]}[1]{\mspace{#1mu}}

\newcommand{\R}[1][n+1]{{\protect\mathbb R}^{#1}}

\newcommand{\Ss}[1][n+1]{{\protect\mathbb S}^{#1}}

\newcommand{\N}{{\protect\mathbb N}}

\newcommand{\eR}{\stackrel{\lower1ex \hbox{\rule{6.5pt}{0.5pt}}}{\msp[3]\R[]}}
\newcommand{\eN}{\stackrel{\lower1ex \hbox{\rule{6.5pt}{0.5pt}}}{\msp[1]\N}}
\newcommand{\eO}{\stackrel{\lower1ex \hbox{\rule{6pt}{0.5pt}}}{\msc O}}

\DeclareMathOperator{\diam}{diam}

\DeclareMathOperator{\graph}{graph}

\newcommand\im{\implies}
\newcommand\ra{\rightarrow}

\newcommand\pa{\partial}
\newcommand\pde[2]{\frac {\partial#1}{\partial#2}}
\newcommand\pd[3]{\frac {\partial#1}{\partial#2^#3}}   %e.g. \pd fxi
\newcommand\pdd[4]{\frac {{\partial\hskip0.15em}^2#1}{\partial {#2^
#3}\,\partial{#2^#4}}}    %e.g. \pdd fxij, Abl. zweiter Ordnung
 
\newcommand\df[2]{\frac {d#1}{d#2}}

\newcommand{\un}{\infty}
\newcommand{\A}{\forall}

\newcommand{\set}[2]{\{\,#1\colon #2\,\}}
\newcommand{\uu}{\cup}
\newcommand{\ii}{\cap}
\newcommand{\uuu}{\bigcup}

\newcommand{\uud}{ \stackrel{\lower 1ex \hbox {.}}{\uu}}
\newcommand{\uuud}[1]{ \stackrel{\lower 1ex \hbox {.}}{\uuu_{#1}}}
\newcommand\su{\subset}
\newcommand\Su{\Subset}

\newcommand{\sminus}[1][28]{\raise 0.#1ex\hbox{$\scriptstyle\setminus$}}

\newcommand{\ol}{\overline}

\newcommand{\wed}{\wedge}

\newcommand{\abs}[1]{\lvert#1\rvert}

\newcommand{\norm}[1]{\lVert#1\rVert}

\newcommand{\spd}[2]{\protect\langle #1,#2\protect\rangle}

\newcommand\ch[3]{\varGamma_{#1#2}^#3}
\newcommand\cha[3]{{\bar\varGamma}_{#1#2}^#3}

\newcommand{\riem}[4]{R_{#1#2#3#4}}
\newcommand{\riema}[4]{{\bar R}_{#1#2#3#4}}

\newcommand{\tit}{\textit}

\newcommand{\tup}{\textup}% text upright

\newcommand{\mc}{\protect\mathcal}
\newcommand{\msc}{\protect\mathscr}

\providecommand{\bysame}{\makebox[3em]{\hrulefill}\thinspace}

\newcommand{\ci}{\cite}

\newcommand{\bt}{\begin{thm}}
\newcommand{\bl}{\begin{lem}}
\newcommand{\bc}{\begin{cor}}
\newcommand{\bd}{\begin{definition}}
\newcommand{\bpp}{\begin{prop}}
\newcommand{\br}{\begin{rem}}
\newcommand{\bn}{\begin{note}}
\newcommand{\be}{\begin{ex}}
\newcommand{\bes}{\begin{exs}}
\newcommand{\bb}{\begin{example}}
\newcommand{\bbs}{\begin{examples}}
\newcommand{\ba}{\begin{axiom}}
\newcommand{\bas}{\begin{assumption}}

\newcommand{\et}{\end{thm}}
\newcommand{\el}{\end{lem}}
\newcommand{\ec}{\end{cor}}
\newcommand{\ed}{\end{definition}}
\newcommand{\epp}{\end{prop}}
\newcommand{\er}{\end{rem}}
\newcommand{\en}{\end{note}}
\newcommand{\ee}{\end{ex}}
\newcommand{\ees}{\end{exs}}
\newcommand{\eb}{\end{example}}
\newcommand{\ebs}{\end{examples}}
\newcommand{\ea}{\end{axiom}}
\newcommand{\eas}{\end{assumption}}

\newcommand{\bp}{\begin{proof}}
\newcommand{\ep}{\end{proof}}
\newcommand{\eps}{\renewcommand{\qed}{}\end{proof}}

\newcommand{\bal}{\begin{align}}

\newcommand{\bi}[1][1.]{\begin{enumerate}[\upshape #1]}
\newcommand{\bia}[1][(1)]{\begin{enumerate}[\upshape #1]}
\newcommand{\bin}[1][1]{\begin{enumerate}[\upshape\bfseries #1]}
\newcommand{\bir}[1][(i)]{\begin{enumerate}[\upshape #1]}
\newcommand{\bic}[1][(i)]{\begin{enumerate}[\upshape\hspace{2\cma}#1]}
\newcommand{\bis}[2][1.]{\begin{enumerate}[\upshape\hspace{#2\parindent}#1]}
\newcommand{\ei}{\end{enumerate}}

\newcommand\ndots{\raise 0.47ex \hbox {,}\hskip0.06em\cdots %
     \raise 0.47ex \hbox {,}\hskip0.06em} 

\newcommand{\q}{\quad}
\newcommand{\qq}{\qquad}

\newcommand{\hp}{\hphantom}

\newcommand\nd{\noindent}

\newskip\Csmallskipamount                                                
\Csmallskipamount=\smallskipamount
\newskip\Cmedskipamount
\Cmedskipamount=\medskipamount
\newskip\Cbigskipamount
\Cbigskipamount=\bigskipamount

\newcommand\cvs{\vspace\Csmallskipamount}   
\newcommand\cvm{\vspace\Cmedskipamount}

\newskip\csa
\csa=\smallskipamount

\newskip\cma
\cma=\medskipamount

\newskip\cba
\cba=\bigskipamount

\newdimen\spt
\newcommand\citem{\cvs\advance\itemno by
1{(\romannumeral\the\itemno})\hskip3pt}
\newcommand{\bitem}{\cvm\nd\advance\itemno by
1{\bf\the\itemno}\hspace{\cma}}

\newcount\itemno
\newcommand{\las}[1]{\label{S:#1}}

\newcommand{\lae}[1]{\label{E:#1}}
\newcommand{\lat}[1]{\label{T:#1}}
\newcommand{\lal}[1]{\label{L:#1}}
\newcommand{\lad}[1]{\label{D:#1}}
\newcommand{\lac}[1]{\label{C:#1}}

\newcommand{\lar}[1]{\label{R:#1}}

\newcommand{\rt}[1]{Theorem~\ref{T:#1}}
\newcommand{\rl}[1]{Lemma~\ref{L:#1}}

\newcommand{\rc}[1]{Corollary~\ref{C:#1}}

\newcommand{\rr}[1]{Remark~\ref{R:#1}}

\newcommand{\re}[1]{\eqref{E:#1}}
\newcommand{\frc}[1]{Corollary~\ref{C:#1} on page~\tup{\pageref{C:#1}}}
\newcommand{\frt}[1]{Theorem~\ref{T:#1} on page~\tup{\pageref{T:#1}}}
\newcommand{\frl}[1]{Lemma~\ref{L:#1} on page~\tup{\pageref{L:#1}}}

\newcommand{\frd}[1]{Definition~\ref{D:#1} on page~\tup{\pageref{D:#1}}}
\newcommand{\fre}[1]{\eqref{E:#1} on page~\tup{\pageref{E:#1}}}

\newcommand{\frs}[1]{Section~\ref{S:#1} on page~\tup{\pageref{S:#1}}} 
\newcommand{\ind}[1]{#1\index{#1}}

\newskip\thmskip
\thmskip=\parindent

\newskip\hsk
\newenvironment{hinw}{\labelsep=0pt\begin{list}{}{\labelsep=0pt\itemindent=0pt\labelwidth=0pt\leftmargin=\parindent\rightmargin=0pt\partopsep=\cba}%
\item\it\nopagebreak\nopagebreak}%
{\end{list}}

\newcommand\bh{\begin{hinw}}
\newcommand{\eh}{\end{hinw}}

\newtheoremstyle{normal}% name
  {\cba}%      Space above, empty = `usual value'
  {\cba}%      Space below
  {}% Body font
  {\thmskip}%Indent amount (empty = no indent, \parindent = para indent)
  {\bfseries}% Thm head font
  {.}%        Punctuation after thm head
  {\hsk}%     Space after thm head: " " = normal interword space;
  {}% Thm head spec

\newtheoremstyle{abschnitt}% name
  {\cba}%      Space above, empty = `usual value'
  {\cba}%      Space below
  {}% Body font
  {\thmskip}% Indent amount (empty = no indent, \parindent = para indent)
  {\bfseries}% Thm head font
  {.}%        Punctuation after thm head
  {\hsk}%     Space after thm head: " " = normal interword space;
  {}% Thm head spec

\newtheoremstyle{italic}% name
  {\cba}%      Space above, empty = `usual value'
  {\cba}%      Space below
  {\itshape}% Body font
  {\thmskip}%  Indent amount (empty = no indent, \parindent = para indent)
  {\bfseries}% Thm head font
  {.}%        Punctuation after thm head
  {\hsk}%     Space after thm head: " " = normal interword space;
  {}% Thm head spec

\newtheoremstyle{aufgaben}% name
  {\cba}%      Space above, empty = `usual value'
  {\cba}%      Space below
  {}% Body font
  {}%         Indent amount (empty = no indent, \parindent = para indent)
  {\normalsize\bfseries}% Thm head font
  {.}%        Punctuation after thm head
  {\hsk}%     Space after thm head: " " = normal interword space;
  {}% Thm head spec

\newtheoremstyle{break}% name
  {\cba}%      Space above, empty = `usual value'
  {\cba}%      Space below
  {\itshape}% Body font
  {}%         Indent amount (empty = no indent, \parindent = para indent)
  {\bfseries}% Thm head font
  {.}%        Punctuation after thm head
  {\newline}% Space after thm head: \newline = linebreak
  {}%         Thm head spec

\swapnumbers
\theoremstyle{italic}
\newtheorem{thm}[subsection]{Theorem}
\newtheorem{lem}[subsection]{Lemma}
\newtheorem{prop}[subsection]{Proposition}
\newtheorem{cor}[subsection]{Corollary}

\theoremstyle{normal}
\newtheorem{rem}[subsection]{Remark}
\newtheorem{definition}[subsection]{Definition}
\newtheorem{example}[subsection]{Example}
\newtheorem{examples}[subsection]{Examples}
\newtheorem{ex}[subsection]{Exercise}
\newtheorem{note}[subsection]{}
\newtheorem{axiom}[subsection]{Axiom}
\newtheorem{assumption}[subsection]{Assumption}

\theoremstyle{aufgaben}
\newtheorem{exs}[subsection]{Exercises}

\swapnumbers

\numberwithin{equation}{section}
\numberwithin{figure}{section}

\newenvironment{textequation}[1][0.8]
{\begin{equation}
\begin{aligned}
\begin{minipage}{#1\linewidth}}
{\end{minipage}
\end{aligned}
\end{equation}
\ignorespacesafterend}

\newcommand{\btext}{\begin{textequation}}
\newcommand{\etext}{\end{textequation}}

\def\hinweis{\@startsection{subsection}{2}%
 \z@{0.7\linespacing\@plus 0.5\linespacing}{0.7\linespacing}%
{\normalfont\itshape\indent}}

\newcounter{hours}\newcounter{minutes}
\newcommand{\printtime}{%
\setcounter{hours}{\time/60}%
\setcounter{minutes}{\time-\value{hours}*60}%
\ifthenelse{\value{minutes}<10}{\thehours :0\theminutes}{\thehours:\theminutes}}

\usepackage[german,english]{babel}
\usepackage{graphicx}
\RequirePackage{amsmath}
\RequirePackage{bm}
\RequirePackage{amssymb}
\RequirePackage{upref}
\RequirePackage{amsthm}
\RequirePackage{enumerate}
\RequirePackage{pb-diagram}
\RequirePackage{amsfonts}
\RequirePackage[mathscr]{eucal}
\RequirePackage{verbatim}
\RequirePackage{xr}
\RequirePackage{graphicx}
\usepackage{calc}
\usepackage{xspace}
\usepackage{dsfont}

\makeatletter
\RequirePackage{color}
\newcommand{\ann}[1]{\renewcommand{\@makefnmark}{\mbox{$^{\color{red}{\@thefnmark}}$}}%
\footnote {#1}}
\makeatother

\RequirePackage{upref}
\RequirePackage{amsthm}
\RequirePackage{enumerate}%\begin{enumerate}[(i)]
\usepackage[mathscr]{eucal}
\usepackage{xr-hyper}

\usepackage{calc}

\newlength{\oddsidemarginlength}
\newlength{\topmarginlength}

\newcounter{numberoflines}
\newcounter{tempcc}
\oddsidemargin=\oddsidemarginlength
\evensidemargin=\oddsidemargin
\topmargin=\topmarginlength
\usepackage[colorlinks=true,linkcolor=blue,citecolor=blue,urlcolor=blue]{hyperref}  

\begin{document}

\flushbottom

%\larger[1]
%\frontmatter

\title{Curvature flows in the sphere}

% author one information
\author{Claus Gerhardt}
\address{Ruprecht-Karls-Universit\"at, Institut f\"ur Angewandte Mathematik,
Im Neuenheimer Feld 294, 69120 Heidelberg, Germany}
%\curraddr{}
\email{\href{mailto:gerhardt@math.uni-heidelberg.de}{gerhardt@math.uni-heidelberg.de}}
\urladdr{\href{http://www.math.uni-heidelberg.de/studinfo/gerhardt/}{http://www.math.uni-heidelberg.de/studinfo/gerhardt/}}
\thanks{This work was supported by the DFG}

% author two information
%\author{}
%\address{}
%\curraddr{}
%\email{}
%\thanks{}
%
\date{\today}
%
% at present the "communicated by" line appears only in ERA and PROC
%\commby{}

%\dedicatory{}

\begin{abstract} 
We consider contracting and expanding curvature flows in $\Ss$. When the flow hypersurfaces are strictly convex we establish a relation between the contracting hypersurfaces and the expanding hypersurfaces which is given by the Gau{\ss} map. The contracting hypersurfaces shrink to a point $x_0$ while the expanding hypersurfaces converge to the equator of the hemisphere $\mc H(-x_0)$. After rescaling, by the same scale factor, the rescaled hypersurfaces converge to the unit spheres with centers $x_0$ \resp $-x_0$ exponentially fast in $C^\un(\Ss[n])$.
\end{abstract}

\maketitle

\tableofcontents

\setcounter{section}{0}

\section{Introduction} 
We consider contracting and expanding curvature flows in $\Ss$. When the flow hypersurfaces are strictly convex we establish a relation between the contracting hypersurfaces and the expanding hypersurfaces which is given by the Gau{\ss} map. Consider monotone curvature functions $F$ being defined in the positive cone $\C_+\su \R[n]$ such that
\begin{equation}
F(1,\ldots,1)=1
\end{equation}
and such that both $F$ and its inverse $\tilde F$ are concave. Let $M(t)$ \resp $\tilde M(t)$ be solutions of the flows
\begin{equation}\lae{1.2}
\dot x=-F\nu
\end{equation}
\resp
\begin{equation}\lae{1.3}
\dot x=\tilde F^{-1}\nu,
\end{equation}
where the initial hypersurfaces $M_0$ \resp $\tilde M_0$ are strictly convex and where $\tilde M_0$ is the polar set of $M_0$, then both flows exist on the maximal time interval $[0,T^*)$, the hypersurfaces $\tilde M(t)$ are the polar hypersurfaces of $M(t)$, and vice versa. The contracting hypersurfaces shrink to a point $x_0$ while the expanding hypersurfaces converge to the equator of the hemisphere $\mc H(-x_0)$. After rescaling, by the same scale factor, the rescaled hypersurfaces 
satisfy uniform estimates in the $C^\un$ topology with uniformly positive principal curvatures. When the curvature function $F$ of the contracting flow is strictly concave, see \frd{3.1} for a precise definition,  or when $F=\tfrac1n H$, then the rescaled hypersurfaces of both flows converge to the unit spheres with centers $x_0$ \resp $-x_0$ exponentially fast in $C^\un(\Ss[n])$.

The class of strictly concave curvature functions comprises the appropriate roots $\s_k$, $2\le k\le n$, of the elementary symmetric polynomials, the functions of class $(K)$, and hence the inverses $\tilde \s_k$ of the $\s_k$, $1\le k\le n$. Proofs of these results concerning strictly concave curvature functions are given in \frs{3}. As a byproduct we also obtain a simple proof that the $\s_k$ are concave. 

Here is a more detailed summary of our results. 
\bt
Let $F\in C^\un(\C_+)$ be a symmetric, monotone and homogeneous of degree $1$ curvature function and assume that both $F$ and its inverse $\tilde F$ are concave. Normalize $F$ such that
\begin{equation}
F(1,\ldots,1)=1
\end{equation}
and consider the curvature flows \re{1.2} \resp \re{1.3} with initial smooth and strictly convex  hypersurfaces $M_0$ \resp $\tilde M_0$, where $\tilde M_0$ is the polar of $M_0$. Then both flows exist in the maximal time interval $[0,T^*)$ with finite $T^*$. The respective flow hypersurfaces are polar sets of each other. The contracting flow hypersurfaces shrink to a point $x_0$ while the expanding hypersurfaces converge to the equator of the hemisphere $\mc H(-x_0)$. The contracting flow is compactly contained in the open hemisphere $\mc H(x_0)$ for   $t_\de\le t<T^*$ while the expanding flow is contained in $\mc H(-x_0)$ for all $0\le t<T^*$. 

Introducing geodesic polar coordinate systems with centers in $x_0$ \resp $-x_0$ and writing the flow hypersurfaces as graphs of a function $u$ \resp $u^*$, then, for any $m\in\N$, we have
\begin{equation}
\abs{u}_{m,\Ss[n]}\le c_m\Theta\qq\A\, t\in [t_\de,T^*)
\end{equation}
\resp
\begin{equation}
\abs{\tfrac\pi2-u^*}_{m,\Ss[n]}\le c_m\Theta\qq\A\, t\in [t_\de,T^*),
\end{equation}
where $\Theta(t,T^*)$ is the solution of the flow \re{1.2} with spherical initial hypersurface and same existence interval.

The rescaled functions 
\begin{equation}\lae{1.7}
u\Theta^{-1}
\end{equation}
\resp
\begin{equation}\lae{1.8}
(\tfrac\pi2-u^*)\Theta^{-1}
\end{equation}
are uniformly bounded $C^\un(\Ss[n])$ and the rescaled principal curvatures  are uniformly positive.

When the curvature function $F$, governing the contracting flow, is strictly concave, or when $F=\tfrac1n H$, then the functions in \re{1.7} \resp \re{1.8} converge to the constant function $1$ in $C^\un(\Ss[n])$ exponentially fast.
\et

Contracting curvature flows have first been considered by Huisken for the mean curvature in Euclidean and Riemannian spaces, \cf \cite{huisken:euclidean,huisken:riemannian}. We are adapting his method of proving an exponential decay for the difference of the principal curvatures to the present situation in order to derive our decay estimates for the rescaled hypersurfaces. Tso proved that contracting hypersurfaces by the Gau{\ss} curvature shrinks the hypersurfaces to a point \cite{tso:gauss}, while Chow proved the contraction to a round point in case of the square root of the scalar curvature and the $n$-th root of the Gau{\ss} curvature, \cf \cite{chow:gauss,chow:scalar}. Andrews, \cite{andrews:euclidean,andrews:riemannian},  considered contracting flows for a class of curvature functions in Euclidean and Riemannian spaces and proved convergence to a point,  boundedness of the rescaled hypersurfaces in the $C^\un$ topology and also convergence to a sphere (or spheres in the Riemannian case), though we do not understand his arguments for the convergence of the rescaled hypersurfaces and consider his proofs to be incorrect. 

Expanding flows, or inverse curvature flows, have been considered in Euclidean  and hyperbolic space \cite{cg90,cg:icf-hyperbolic,cg:icfp,scheuer:icfp}. Recently,  inverse curvature flows have been studied in $\Ss$ by Makowski and Scheuer \cite{scheuer:icfs} who proved convergence to a hemisphere in $C^{1,\al}$.

\br
Our results for the contracting flows are also valid in $\R$.
\er
\section{Definitions and notations}
The main objective of this section is to state the equations of Gau{\ss}, Codazzi,
and Weingarten for  hypersurfaces $M$ in a \di {(n+1)} Riemannian
manifold
$N$.  Geometric quantities in $N$ will be denoted by
$(\bar g_{ \al \bet})$, $(\riema  \al \bet \ga \de)$, etc., and those in $M$ by $(g_{ij}), 
(\riem ijkl)$, etc. Greek indices range from $0$ to $n$ and Latin from $1$ to $n$; the
summation convention is always used. Generic coordinate systems in $N$ resp.
$M$ will be denoted by $(x^ \al)$ resp. $(\x^i)$. Covariant differentiation will
simply be indicated by indices, only in case of possible ambiguity they will be
preceded by a semicolon, i.e., for a function $u$ in $N$, $(u_ \al)$ will be the
gradient and
$(u_{ \al \bet})$ the Hessian, but e.g., the covariant derivative of the curvature
tensor will be abbreviated by $\riema  \al \bet \ga{ \de;\e}$. We also point out that
\begin{equation}
\riema  \al \bet \ga{ \de;i}=\riema  \al \bet \ga{ \de;\e}x_i^\e
\end{equation}
with obvious generalizations to other quantities.

Let $M$ be a  $C^2$-hypersurface 
with  normal $\n$.

In local coordinates, $(x^ \al)$ and $(\x^i)$, the geometric quantities of the
hypersurface $M$ are connected through the following equations
\begin{equation}\lae{2.2}
x_{ij}^ \al= -h_{ij}\n^ \al
\end{equation}
the so-called \tit{Gau{\ss} formula}. Here, and also in the sequel, a covariant
derivative is always a \tit{full} tensor, i.e.,
\begin{equation}
x_{ij}^ \al=x_{,ij}^ \al-\ch ijk x_k^ \al+ \cha  \bet \ga \al x_i^ \bet x_j^ \ga.
\end{equation}
The comma indicates ordinary partial derivatives.

In this implicit definition the \tit{second fundamental form} $(h_{ij})$ is taken
with respect to $-\n$.

The second equation is the \tit{Weingarten equation}
\begin{equation}
\n_i^ \al=h_i^k x_k^ \al,
\end{equation}
where we remember that $\n_i^ \al$ is a full tensor.

Finally, we have the \tit{Codazzi equation}
\begin{equation}
h_{ij;k}-h_{ik;j}=\riema \al \bet \ga \de\n^ \al x_i^ \bet x_j^ \ga x_k^ \de
\end{equation}
and the \tit{Gau{\ss} equation}
\begin{equation}
\riem ijkl= \{h_{ik}h_{jl}-h_{il}h_{jk}\} + \riema  \al \bet\ga \de x_i^ \al x_j^ \bet
x_k^ \ga x_l^ \de.
\end{equation}

When we consider hypersurfaces $M\su S^{n+1}$ to be embedded in $\R[n+2]$, we label the coordinates in $\R[n+2]$ as $(x^a)$, i.e., indices $a,b,c,...$ always run through $n+2$ values either from $1$ to $n+2$ or from $0$ to $n+1$. 

At the end of this section let us state some evolution equations satisfied by solutions of the curvature flows
\begin{equation}\lae{2.7.314} 
\dot x=-\F\nu
\end{equation}
in a Riemannian space form $N=N^{n+1}$ with curvature $K_N$. Here $\F=\F(F)$.
\bl\lal{1.2.314}
The term $\F$ evolves according to the equation
\begin{equation}\lae{1.2.314}
\begin{aligned}
{\F}^\prime-\dot\F F^{ij}\F_{ij}&=
\dot \F
F^{ij}h_{ik}h_j^k \F\\
&\hp{=}\;+ K_N\dot\F F^{ij}g_{ij}\F,
\end{aligned}
\end{equation}where
\begin{equation}
\F^{\prime}=\frac{d}{dt}\F
\end{equation}
and
\begin{equation}
\dot\F=\frac{d}{dr}\F(r).
\end{equation}
\el
For a proof see \cite[Lemma 2.3.4]{cg:cp}.

Assume that the flow hypersurfaces are written as graphs in a geodesic polar coordinate system. Define $v$ by
\begin{equation}
v^{-1}=\spd{\pde{}{x^0}}\nu
\end{equation}
and let $\h=\h(r)$ be a positive solution of the equation
\begin{equation}
\dot\h=-\frac{\bar H}n\h,
\end{equation}
where $\bar H$ is the mean curvature of the slices $\{x^0=r\}$, then
\begin{equation}
\chi=v\h(u)
\end{equation}
satisfies the equation
\begin{equation}\lae{1.5.314}
\dot\chi-\dot\F F^{ij}\chi_{ij}=-\dot\F F^{ij}h_{ik}h^k_j\chi-2\chi^{-1}\dot\F F^{ij}\chi_i\chi_j+\{\dot\F F+\F\}\frac{\bar H}nv\chi,
\end{equation}
\cf \cite[Lemma 5.8]{cg:spaceform} 
\bl\lal{1.3.314}
Let $N$ be a space of constant curvature $K_N$, then the second fundamental form of the curvature flow \re{2.7.314} satisfies the  parabolic equations
\begin{equation}\lae{1.6.314}
\begin{aligned}
\dot h_i^j-\dot\F F^{kl}h_{i;kl}^j&=\dot\F F^{kl}h_{rk}h_l^rh_i^j-\dot\F F
h_{ri}h^{rj}+\F h_i^kh_k^j\\
&\hp{=}\;+\Ddot\F F_iF^j+\dot\F
F^{kl,rs}h_{kl;i}h_{rs;}^{\hphantom{rs;}j}\\
&\hp{=}\;+K_N\{\F\de^j_i+\dot\F F\de^j_i-\dot\F F^{kl}g_{kl}h^j_i\}.
\end{aligned}
\end{equation}
and 
\begin{equation}\lae{1.6.1.314}
\begin{aligned}
\dot h_{ij}-\dot\F F^{kl}h_{ij;kl}&=\dot\F F^{kl}h_{rk}h_l^rh_{ij}-\dot\F F
h_{ri}h^{r}_j-\F h_i^kh_{kj}\\
&\hp{=}\;+\Ddot\F F_iF_j+\dot\F
F^{kl,rs}h_{kl;i}h_{rs;j}\\
&\hp{=}\;+K_N\{\F g_{ij}+\dot\F Fg_{ij}-\dot\F F^{kl}g_{kl}h_{ij}\}.
\end{aligned}
\end{equation}
\el
For a proof see \cite[Lemma 2.4.3]{cg:cp}.

\bl\lal{1.4.314}
Let $h_{ij}$ be invertible and set $(\tilde h^{ij})=(h_{ij})^{-1}$, then the mixed tensor $\tilde h^i_j$ satisfies the evolution equation
\begin{equation}\lae{1.7.314}
\begin{aligned}
\dot{\tilde h}^i_j-\dot\F F^{kl}\tilde h^i_{j;kl}&=\\
&\msp[-120]-\;\dot\F F^{kl}h_{kr}h^r_l\tilde h^i_j+\{\dot\F F -(\F-\tilde f)\}\de^i_j\\
&\msp[-120]-K_N\{\dot\F F+\F \}\tilde h_{kj}\tilde h^{ki} +K_N\dot\F F^{kl}g_{kl}\tilde h^i_j\\
&\msp[-120]-\{\dot\F F^{pq,kl}h_{pq;r}h_{kl;s}+2\dot\F F^{kl}\tilde h^{pq}h_{pk;r}h_{ql;s}+\Ddot\F F_rF_s\}\tilde h^{is}\tilde h^r_j.
\end{aligned}
\end{equation} 
\el

\section{Curvature functions}\las{3}

\bd\lad{3.1}
Let $F\in C^2(\C)$ be a symmetric, homogeneous of degree $1$, monotone and concave curvature function. We call $F$ \tit{strictly concave}, if the multiplicity of the eigenvalue $\lam=0$ for $D^2F(\ka)$ is one for all $\ka\in\C$.
\ed
We shall show that the $k$-th root of the elementary symmetric polynomials $H_k$, $2\le k\le n$, are strictly concave. This will also offer a simple independent proof of the concavity of the $k$-th root of $H_k$.

The $H_k$ are defined in the connected component $\C_k$ of the cone
\begin{equation}
\{H_k>0\}
\end{equation}
containing $\C_+$. The cones are monotonely ordered
\begin{equation}
\C_+=\C_n\su\cdots\su\C_1,
\end{equation}
\cf \cite[Section 2]{huisken-s:cone}.

\bt
The curvature functions
\begin{equation}
\s_k=H_k^\frac1{k},\qq 2\le k\le n,
\end{equation}
are strictly concave.
\et
\bp
The proof relies on the concavity of the functions
\begin{equation}
Q_k=\frac{H_{k+1}}{H_k},\qq 1\le k\le n-1.
\end{equation}
A proof of this fact can be found in \cite[Theorem 2.5]{huisken-s:cone}. There, it also proved that the $Q_k$ are strictly concave in $\C_+$.

For the proof of the theorem we shall use induction with respect to $k$. A proof that $\s_2$ is strictly concave is given in the lemma below.

Thus, let us assume that $\s_k$, $2\le k<n$, is already strictly concave. Define
\begin{equation}
F=\s_{k+1},
\end{equation}
then
\begin{equation}
\begin{aligned}
F_{ij}&= (\tfrac1{k+1}-1)\tfrac1{k+1}H^{\frac1{k+1}-2}_{k+1} H_{k+1,i}H_{k+1,j}\\
&\hp{=}\;+\tfrac1{k+1}H^{\frac1{k+1}-1}_{k+1}H_{k+1,ij}.
\end{aligned}
\end{equation}
Here, the indices denote partial derivatives. Then the concavity of $F$ is equivalent to the relation
\begin{equation}\lae{3.7.1}
H_{k+1,ij}\le (1-\tfrac1{k+1})H^{-1}_{k+1}H_{k+1,i}H_{k+1,j}.
\end{equation}
We shall prove this inequality by induction and also
\begin{equation}\lae{3.8.1}
H_{k+1,ij}\xi^i\xi^j< (1-\tfrac1{k+1})H^{-1}_{k+1}H_{k+1,i}\xi^iH_{k+1,j}\xi^j\qq\A\, \ka\not\sim\xi\in\R[n],
\end{equation}
where $\xi\not=0$ and where $\ka\sim\xi$ means that
\begin{equation}
\xi=\lam \ka.
\end{equation}

Let $\f$ be defined by
\begin{equation}
\f=Q_{k},
\end{equation}
then
\begin{equation}
\begin{aligned}
H_{k+1,ij}=\f_{ij}H_k+\f_iH_{k,j}+\f_jH_{k,i}+\f H_{k,ij}.
\end{aligned}
\end{equation}
The argument $\ka\in\C$ is obviously an eigenvector of $D^2F(\ka)$ with eigenvalue $0$. Hence, let $\ka\not\sim\xi\in\R[n]$ be arbitrary, $\xi\not=0$, then we deduce
\begin{equation}
\begin{aligned}
H_{k+1,ij}\xi^i\xi^j&=\f_{ij}\xi^i\xi^jH_k+2\f_i\xi^iH_{k,j}\xi^j+\f H_{k,ij}\xi^i\xi^j\\
&<2\f_i\xi^iH_{k,j}\xi^j+\frac{H_{k+1}}{H_k}(1-\tfrac1k)H_k^{-1}(H_{k,i}\xi^i)^2,
\end{aligned}
\end{equation}
where we used the concavity of $\f$ and the assumption \re{3.8.1} for the function $H_k$.

From the relation
\begin{equation}
H_{k+1,i}\xi^i=\f_i\xi^iH_k+\f H_{k,i}\xi^i
\end{equation}
we obtain
\begin{equation}
\begin{aligned}
\f_i\xi^iH_{k,j}\xi^j=H_k^{-1}H_{k+1,i}\xi^iH_{k,j}\xi^j-\frac{H_{k+1}}{H_k^2}(H_{k,i}\xi^i)^2
\end{aligned}
\end{equation}
yielding
\begin{equation}
\begin{aligned}
H_{k+1,ij}\xi^i\xi^j&<2H_k^{-1}H_{k+1,i}\xi^iH_{k,j}\xi^j-2\frac{H_{k+1}}{H_k^2}(H_{k,i}\xi^i)^2\\
&\hp{<}\;+\frac{k-1}k\frac{H_{k+1}}{H_k^2}(H_{k,i}\xi^i)^2\\
&\le \tfrac k{k+1}H_{k+1}^{-1}(H_{k+1,i}\xi^i)^2 +\frac{k+1}k\frac{H_{k+1}}{H_k^2}(H_{k,i}\xi^i)^2\\
&\hp{\e}\; -2\frac{H_{k+1}}{H_k^2}(H_{k,i}\xi^i)^2+\frac {k-1}k \frac{H_{k+1}}{H_k^2}(H_{k,i}\xi^i)^2\\
&=(1-\tfrac 1{k+1})H_{k+1}^{-1}(H_{k+1,i}\xi^i)^2 .
\end{aligned}
\end{equation}
The lemma below will complete the proof of the theorem.
\ep
\bl
The function $\s_2$ is strictly concave.
\el
\bp
We shall first prove that $F=\s_2$ is concave. We use the same technique as in the proof of the theorem above and shall verify that the inequality \re{3.7.1} is satisfied for $F$. Define
\begin{equation}
\f=\frac{H_2}H
\end{equation}
and let $\xi\in\R[n]$, then
\begin{equation}
\begin{aligned}
H_{2,ij}\xi^i\xi^j&\le 2H^{-1}H_{2,i}\xi^iH_j\x^j-2\frac{H_2}{H^2}(H_i\xi^i)^2\\
&\le \tfrac12H_2^{-1}(H_{2,i}\xi^i)^2+2\frac{H_2}{H^2}(H_i\xi^i)^2-2\frac{H_2}{H^2}(H_i\xi^i)^2\\
&=\tfrac12H_2^{-1}(H_{2,i}\xi^i)^2,
\end{aligned}
\end{equation}
hence $F=\s_2$ is concave.

To prove that $F$ is strictly concave, assume there exists $0\not=\xi\in\R[n]$ such that
\begin{equation}
F_{ij}(\ka)\xi^j=0\q\wed\q \ka_i\xi^i=0.
\end{equation}
For simplicity let us define
\begin{equation}
F=\sqrt{H^2-\abs A^2},
\end{equation}
then
\begin{equation}
F_i=F^{-1}(H-\ka_i)
\end{equation}
and
\begin{equation}
F_{ij}=-F^{-3}(H-\ka_i)(H-\ka_j)+F^{-1}(1-\de_{ij}).
\end{equation}
Define $\s$ by
\begin{equation}
\s=\sum_i\xi^i,
\end{equation}
then
\begin{equation}\lae{3.22.1}
(\s-\xi^i)F^2=(H-\ka_i)H\s.
\end{equation}
Summing over $i$ yields
\begin{equation}
(n-1)\s F^2=(n-1)H^2\s,
\end{equation}
and hence we deduce 
\begin{equation}
\s=0
\end{equation}
for otherwise we get a contradiction. But when $\s=0$, we infer from \re{3.22.1}
\begin{equation}
\xi^i F^2=0,
\end{equation}
a contradiction.
\ep
Now, we want to prove that the inverses $\tilde \s_k$ of $\s_k$, $1\le k\le n$, are also strictly concave. This will follow from the fact that they are of class $(K)$.
\bd
A symmetric curvature function $F\in C^{2,\alpha}(\C_+)\ii C^0(\bar \C_+)$,
positively homogeneous of degree $d_0>0$, is said to be of class $(K)$, if
\begin{equation}\lae{1.2.1.1}
F_i=\pd F{\kappa}i>0\q \text{in } \C_+,
\end{equation}
which is also referred to as $F$ to be \tit{\ind{strictly monotone}},
\begin{equation}\lae{1.2.28}
\fmo{\pa \C_+}=0,
\end{equation}
and
\begin{equation}\lae{1.2.29}
F^{ij,kl}\h_{ij}\h_{kl}\le F^{-1}(F^{ij}\h_{ij})^2-F^{ik}\tilde
h^{jl}\h_{ij}\h_{kl}\qq\A\,\h\in\msc S,
\end{equation}
or, equivalently, if we set $\hat F=\log F$,
\begin{equation}\lae{1.2.30}
\hat F^{ij,kl}\h_{ij}\h_{kl}\le -\hat F^{ik}\tilde
h^{jl}\h_{ij}\h_{kl}\qq\A\,\h\in\msc S,
\end{equation}
where $F$ is evaluated at $(h_{ij})$ and $(\tilde h^{ij})$ is the inverse of $(h_{ij})$.
\ed
Note that we only consider curvature functions which are homogeneous of degree $1$.
\br
 The inverses $\tilde\s_k$ of $\s_k$, $1\le k\le n$, are of class $(K)$, \cf \cite[Chapter 2.2]{cg:cp}, especially Lemma 2.2.11.
\er
\bl
Let $F\in (K)$ be homogenous of degree $1$, then $F$ is strictly concave. 
\el
\bp
The Hessian of $F$ satisfies the inequality
\begin{equation}
\pdd F{\ka}ij\le F^{-1} F_iF_j-F_i\ka_j^{-1}\de_{ij},
\end{equation}
\cf \cite[inequality (2.2.9)]{cg:cp}. The right-hand side is strictly negative definite unless evaluated for a multiple of $\ka$. Indeed, let $\ka\not\sim\xi\in\R[n]$, $\xi\not=0$,  then, using Schwarz's inequality, we deduce
\begin{equation}
\begin{aligned}
F_i\xi^i&=\sum_i F_i^{\frac12}\ka_i^{\frac12}F_i^{\frac12}\ka_i^{-\frac12}\xi^i\\
&\le \big(\sum_iF_i\ka_i\big)^{\frac12} \big(\sum_iF_i\ka_i^{-1}\abs{\xi^i}^2\big)^{\frac12}=F^{\frac12} \big(\sum_iF_i\ka_i^{-1}\abs{\xi^i}^2\big)^{\frac12},
\end{aligned}
\end{equation}
where the inequality is a strict inequality unless
\begin{equation}
\ka_i^{-\frac12}\xi^i=\lam \ka_i^\frac12\qq\A\, i,
\end{equation}
or equivalently,
\begin{equation}
\xi^i=\lam\ka_i\qq\A\, i.
\end{equation}
\ep

\section{Polar sets and dual flows}

Let $M\su S^{n+1}$ be a connected, closed, immersed, strictly convex hypersurface given by an immersion
\begin{equation}
x:M_0\ra M\su S^{n+1},
\end{equation}
then $M$ is embedded, homeomorphic to $S^n$, contained in an open hemisphere and is the boundary of a convex body $\hat M\su S^{n+1}$, \cf \ci{docarmo}.

Considering $M$ as a codimension $2$ submanifold of $\R[n+2]$ such that
\begin{equation}\lae{2.2b}
x_{ij}=-g_{ij}x-h_{ij}\tilde x,
\end{equation}
where $\tilde x\in T_x(\R[n+2])$ represents the exterior normal vector $\nu\in T_x(S^{n+1})$, we   proved in \cite[Theorem 9.2.5]{cg:cp} that the mapping
\begin{equation}\lae{2.3b}
\tilde x:M_0\ra S^{n+1}
\end{equation}
is an embedding of a strictly convex, closed, connected hypersurface $\tilde M$. We called this mapping the \tit{Gau{\ss} map} of $M$. More precisely, we proved

\bt\lat{3.1}
Let $x:M_0\ra M\su S^{n+1}$ be a closed, connected, strictly convex hypersurface of class $C^m$, $m\ge 3$, then the Gau{\ss} map $\tilde x$ in \re{2.3b} is the embedding of a closed, connected, strictly convex hypersurface $\tilde M\su S^{n+1}$ of class $C^{m-1}$. 

Viewing $\tilde M$ as a codimension $2$ submanifold in $\R[n+2]$, its Gaussian formula is
\begin{equation}\lae{2.25b}
\tilde x_{ij}=-\tilde g_{ij}\tilde x-\tilde h_{ij} x,
\end{equation}
where $\tilde g_{ij}$, $\tilde h_{ij}$ are the metric and second fundamental form of the hypersurface $\tilde M\su  S^{n+1}$, and $x=x(\xi)$ is the embedding of $M$ which also represents the exterior normal vector of $\tilde M$. The second fundamental form $\tilde h_{ij}$ is defined with respect to the interior normal vector. 

The second fundamental forms of $M$, $\tilde M$ and the corresponding principal curvatures $\ka_i$, $\tilde \ka_i$ satisfy
\begin{equation}\lae{2.26b}
h_{ij}=\tilde h_{ij}=\spd{\tilde x_i}{x_j}
\end{equation}
and
\begin{equation}\lae{2.27b}
\tilde \ka_i=\ka_i^{-1}.
\end{equation}
\et

If $M$ is supposed to satisfy a curvature equation of the form
\begin{equation}\lae{3.7}
\fv FM=f(\nu),
\end{equation}
where $F$ is a curvature function defined in $\C_+$, $F=F(\ka_i)$, $F$ symmetric, monotone,  homogenous of degree $1$ and smooth (for simplicity), $F\in C^\un(\C_+)$, then the polar set $\tilde M$ of $M$ satisfies the equation
\begin{equation}\lae{3.9}
\fv {\tilde F}{\tilde M}=\frac 1{f(x)},
\end{equation}
where $\tilde F$ is the inverse of $F$,
\begin{equation}
\tilde F(\ka_i)=\frac1{F(\ka_i^{-1})}.
\end{equation}

One may consider the equation \re{3.7} and \re{3.9} to describe dual problems. This duality is also valid in case of curvature flows.

Let $x=x(t,\xi)$ be a solution of the curvature flow
\begin{equation}\lae{3.10}
\dot x=-\F\nu,
\end{equation}
where $\F=\F(r)$ is a smooth real, strictly monotone function defined on $\R[]_+$ and where the $F$ on the right-hand side of \re{3.10} is an abbreviation for
\begin{equation}
\F=\F(F).
\end{equation}
Assume that the flow in \re{3.10} with initial strictly convex hypersurface $M_0$ exists on a maximal time interval $[0,T^*)$ and that the flow hypersurfaces $M(t)$ are strictly convex. Let us consider the flow as flow in $\R[n+2]$, then \re{3.10} takes the form
\begin{equation}
\dot x=-\F\tilde x,
\end{equation}
since
\begin{equation}
\spd{\dot x}x=0,
\end{equation}
$\tilde x$ represents $\nu$ in $T_x(\R[n+2])$ and
\begin{equation}
T_x(\R[n+2])=T_x(\Ss)\oplus \langle x\rangle.
\end{equation}

We also note that $x$ is the normal to $\tilde M$ and that the Weingarten equation has the form
\begin{equation}
x_j=\tilde h^k_j\tilde x_k,
\end{equation}
\cf \cite[Lemma 9.2.4]{cg:cp}. Furthermore, we have, \cf \cite[equ. (9.2.36)]{cg:cp},
\begin{equation}
\spd x{\tilde x}=0,
\end{equation}
and we infer
\begin{equation}\lae{3.17}
\spd x{\dot{\tilde x}}=\F,
\end{equation}
\begin{equation}
\spd{x_j}{\tilde x}=0,
\end{equation}
\begin{equation}
0=\spd{\dot x_j}{\tilde x}+\spd{x_j}{\dot{\tilde x}},
\end{equation}
as well as
\begin{equation}
\dot x_j=-\F_j\tilde x-\F \tilde x_j
\end{equation}
in view of \re{3.10}. Thus, we deduce
\begin{equation}\lae{3.21}
\spd{\dot{\tilde x}}{x_j}=-\spd{\dot x_j}{\tilde x}=\F_j.
\end{equation}

Taking \re{3.17}, \re{3.21} and
\begin{equation}
\spd{\dot{\tilde x}}{\tilde x}=0
\end{equation}
into account we finally conclude
\begin{equation}
\begin{aligned}
\dot{\tilde x}&=\F x +\F^mx_m\\
&=\F x+ \F^m\tilde h^k_m\tilde x_k,
\end{aligned}
\end{equation}
where 
\begin{equation}
\F^m=g^{mj}\F_j.
\end{equation}
The corresponding flow equation in $\Ss$ has the form
\begin{equation}\lae{3.25}
\dot{\tilde x}=\F\tilde \nu + \F^m\tilde h^k_m\tilde x_k.
\end{equation}
Let $t_0\in[0,T^*)$ and introduce polar coordinates with center in the convex body defined by $\tilde M(t_0)$, then, for $t_0\le t<t_0+\e$, $\tilde M(t)$ can be written as graph over $\Ss[n]$
\begin{equation}
\tilde M(t)=\graph \tilde u_{|_{\Ss[n]}},
\end{equation}
and we obtain the scalar curvature flow equation
\begin{equation}
\dot{\tilde u}=\frac{d\tilde u}{d t}=\F\tilde v^{-1}+\F^m\tilde h^k_m\tilde u_k
\end{equation}
by looking at the $0$-th component of \re{3.25}, where
\begin{equation}
\begin{aligned}
\tilde v^2&=1+\frac1{\sin^2\tilde u}\s^{ij}\tilde u_i\tilde u_j\\
&\equiv 1+\abs{D\tilde u}^2
\end{aligned}
\end{equation}
and 
\begin{equation}
\tilde \nu=\tilde v^{-1}(1,-\check{\tilde u}^i)
\end{equation}
such that
\begin{equation}
\abs{D\tilde u}^2=\check{\tilde u}^i\tilde u_i.
\end{equation}

The partial derivative of $\tilde u$ with respect to $t$ then satisfies
\begin{equation}\lae{3.31}
\begin{aligned}
\pde {\tilde u}t&=\dot{\tilde u}-\tilde u_i\dot{\tilde x}^i\\
&=\F\tilde v^{-1}+\F^m\tilde h^k_m\tilde u_k+\F\tilde v^{-1}\abs{D\tilde u}^2-\F^m\tilde h^k_m\de^i_k\tilde u_i\\
&=\tilde v\F.
\end{aligned}
\end{equation}
This is exactly the scalar curvature equation, by considering the partial derivative of $\tilde u$ with respect to $t$, of the  flow equation
\begin{equation}
\dot{\tilde x}=\F\tilde \nu,
\end{equation}
where
\begin{equation}
\F=\F(F)=\F(\tilde F^{-1}),
\end{equation}
$\tilde F$ is the inverse of $F$, i.e., when the $M(t)$ satisfy the inverse curvature flow equation
\begin{equation}
\dot x=\frac1F\nu
\end{equation}
then the polar sets $\tilde M(t)$ satisfy the direct flow equation
\begin{equation}
\dot{\tilde x}=-\tilde F\tilde\nu
\end{equation}
and vice versa.
\bt
Let $\F\in C^\un(\R[]_+)$ be strictly monotone, $\dot\F>0$, and let $F\in C^\un(\C_+)$ be a symmetric, monotone, homogeneous of degree $1$ curvature function such that
\begin{equation}
\fv F{\C_+}>0
\end{equation}
and such that the flows
\begin{equation}
\dot x=-\F(F)\nu
\end{equation}
\resp
\begin{equation}\lae{3.38}
\dot{\tilde x}=\F(\tilde F^{-1})\tilde\nu
\end{equation}
with initial strictly convex hypersurfaces $M_0$ \resp $\tilde M_0$ exist on maximal time intervals $[0,T^*)$ \resp $[0,\tilde T^*)$, where the flow hypersurfaces are strictly convex. Let $M(t)$ \resp $\hat M(t)$  be the corrsponding flow hypersurfaces then $T^*=\tilde T^*$ and $\hat M(t)=\tilde M(t)$.
\et
\bp
In view of the symmetry involved it suffices to prove
\begin{equation}
T^*\le \tilde T^*\q\wed\q \hat M(t)=\tilde M(t)\q \A\, t\in [0,T^*) .
\end{equation}
Let $\Lam$ be defined by
\begin{equation}
\Lam=\set{T\in [0,T^*)}{\tilde M(t) \tup{ solves \re{3.38} }\A\, t\in [0,T]}.
\end{equation}
$\Lam$ is evidently not empty, since a small one-sided neighbourhood of $0$ belongs to $\Lam$ in view of the uniqueness of the solution of the scalar curvature flow
\begin{equation}
\pde{\tilde u}t=\tilde v\F
\end{equation}
with given initial value and the arguments leading to \re{3.31}.

By the same reasoning $\Lam$ is obviously open, while the closedness of $\Lam$ is trivial.
\ep

We shall employ this duality by choosing
\begin{equation}
\F(r)=-r^{-1},
\end{equation}
i.e., we shall study and solve inverse curvature flows and direct curvature flows simultaneously using their specific properties to our advantage.

\section{First estimates}
From now on we assume that both $F$, $\tilde F$ are concave and that
\begin{equation}
F(1,\ldots,1)=\tilde F(1,\ldots,1)=1.
\end{equation}
$\F$ is defined by
\begin{equation}
\F(r)=-r^{-1}
\end{equation}
and we consider the curvature flows
\begin{equation}\lae{4.3}
\dot x=-\F\nu
\end{equation}
and the dual flow
\begin{equation}\lae{4.4}
\dot x=-\tilde F\nu
\end{equation}
with initial hypersurfaces $M_0$ \resp $\tilde M_0$. Both flows exist on a maximal time interval $[0,T^*)$. Let us start with some important estimates.
\bl\lal{4.1}
Let $M(t)$ be a solution of the flow \re{4.3}, then the principal curvatures are uniformly bounded during the evolution
\begin{equation}
\ka_i\le \const.
\end{equation}
\el
\bp
Label the $\ka_i$ such that
\begin{equation}
\ka_1\le\cdots\le \ka_n.
\end{equation}
Then we can pretend that
\begin{equation}
\ka_n=h^n_n
\end{equation}
is smooth and that we apply the parabolic maximum principle to $h^n_n$ in equation \fre{1.6.314}, for details see the proof of \cite[Lemma 3.3.3]{cg:cp}.

Thus, fix $0<T<T^*$ and let $(t_0,\xi_0)$, $0<t_0\le T$,  be a point such that
\begin{equation}
h^n_n(t_0,\xi_0)=\sup_{t\in [0,T]}\sup_{M(t)}h^n_n.
\end{equation}
Then we deduce from \re{1.6.314}
\begin{equation}
\begin{aligned}
0&\le \dot\F F^{kl}h_{ki}h^i_l\msp[1]h^n_n-\dot\F F\abs{h^n_n}^2-F^{-1}\abs{h^n_n}^2-K_N\dot\F F^{ij}g_{ij}h^n_n\\
&\le -F^{-1}\abs{h^n_n}^2-K_N\dot\F F^{ij}g_{ij}h^n_n,
\end{aligned}
\end{equation}
a contradiction, i.e., the maximum is attained at $t=0$.
\ep

\bl
Let $\tilde M(t)$ be a solution of the flow \re{4.4}, then there exists $0<\e_0<\frac1n$ such that
\begin{equation}
\e_0\tilde\ka_n\le \e_0\tilde H\le \tilde \ka_1
\end{equation}
during the evolution, where the principal curvature are labelled
\begin{equation}
\tilde\ka_1\le\ldots\tilde\ka_n
\end{equation}
and where
\begin{equation}
\tilde H=\sum_i\tilde\ka_i.
\end{equation}
\el
\bp
We apply a maximum principle for tensors which was originally proved by Hamilton \cite[Theorem 9.1]{hamilton:ricci} and later generalized by Andrews \cite[Theorem 3.2]{andrews:pinching}. Looking at the equation \fre{1.6.1.314} we deduce that the tensor
\begin{equation}
T_{ij}=\tilde h_{ij}-\e_0\tilde H\tilde g_{ij}
\end{equation}
satisfies the equation 
\begin{equation}
\begin{aligned}
\dot T_{ij}-\tilde F^{kl}T_{ij;kl}&=\tilde F^{kl}\tilde h_{kr}\tilde h^r_lT_{ij}-2\tilde F\tilde h^k_i\tilde h_{kj}+2\e_0\tilde F\tilde H\tilde h_{ij}\\
&\q+2K_N\tilde F(1-\e_0n)\tilde g_{ij}-K_N\tilde F^{kl}\tilde g_{kl}T_{ij}\\
&\q+\tilde F^{kl,rs}\tilde h_{kl;i}\tilde h_{rs;j}-\e_0\tilde F^{kl,rs}\tilde h_{kl;i}\tilde h_{rs;j}\tilde g^{ij}\\
&\equiv N_{ij}+\tilde N_{ij},
\end{aligned}
\end{equation}
where
\begin{equation}
\tilde N_{ij}=\tilde F^{kl,rs}\tilde h_{kl;i}\tilde h_{rs;j}-\e_0\tilde F^{kl,rs}\tilde h_{kl;p}\tilde h_{rs;q}\tilde g^{pq}\tilde g_{ij}. 
\end{equation}

Hamilton's maximum principle then has the form: if the tensor $T_{ij}$ is strictly positive definite at time $t=0$ and if the right-hand side satisfies the so-called null eigenvector condition, i.e.,  $T_{ij}\ge0$ and $T_{ij}\h^j=0$ implies
\begin{equation}
N_{ij}\h^i\h^j+\tilde N_{ij}\h^i\h^j\ge 0,
\end{equation}
then $T_{ij}>0$ during the evolution.

However, the term $\tilde N_{ij}$ does not satisfy a null eigenvector condition in general. Andrews therefore proved in \cite[Theorem 3.2]{andrews:pinching} that the conclusion is still valid if $\tilde N_{ij}$ satisfies the weaker condition
\begin{equation}\lae{4.17}
\tilde N_{ij}\h^i\h^j+\sup_{\C=(\C^r_k)}2\tilde F^{kl}(2\C^r_lT_{ir;k}\h^i-\C^r_k\C^s_lT_{rs})\ge0.
\end{equation}

Moreover, he proved that the weaker condition is satisfied by the present tensor $\tilde N_{ij}$, \cf \cite[Theorem 4.1]{andrews:pinching}, provided $\tilde F$ and $F$ are both concave, \cf \cite[Corollary 2.4]{andrews:pinching}. Hence, the maximum principle can be applied provided $N_{ij}$ satisfies the null eigenvector condition, which can be easily verified by choosing coordinates such that
\begin{equation}
\tilde g_{ij}=\de_{ij}\q\wed\q \h^i=\de^i_1,
\end{equation}
and using the fact that $K_N\ge0$. Of course $\e_0$ has to be sufficiently small such that $T_{ij}>0$ at time $t=0$.
\ep
\section{Contracting flows: Convergence to a point}\las{5}

From now on we are mainly considering contracting flows. To facilitate notation we drop any tildes, i.e., the curvature function involved is denoted by $F$ and the flow equation is
\begin{equation}\lae{5.1}
\dot x=-F\nu.
\end{equation}
In view of the results in the previous section there exist uniform positive constants $c_1$ and $c_2$ such that the principal curvatures
\begin{equation}
\ka_1\le \cdots\le \ka_n
\end{equation}
satisfy the estimates
\begin{equation}
c_1\le \ka_1
\end{equation}
and
\begin{equation}\lae{5.4}
\ka_n\le c_2\ka_1.
\end{equation}

When the initial  hypersurface is a geodesic sphere the flow hypersurfaces are all spheres with the same center and their radii $\Theta=\Theta(t)$ satisfy the equation
\begin{equation}\lae{5.5}
\dot\Theta=-\frac{\cos\Theta}{\sin\Theta}.
\end{equation}
The spherical flows exist only for a finite time, hence the flow \re{5.1} exists only for a finite time and there exists a spherical flow  $\Theta=\Theta(t,T^*)$  which shrinks to a point when $t$ approaches $T^*$, where $T^*$ is the maximal existence time for the flow \re{5.1}. These claims can be immediately deduced by looking at initial spheres $M_1$ \resp $M_2$ such that the initial convex body $\hat M_0$, where $M_0$ is the initial hypersurface of the general flow, satisfies
\begin{equation}\lae{5.6}
B_1\Su\hat M_0\Su B_2,
\end{equation}
where
\begin{equation}
\pa B_i=M_i,\qq i=1,2.
\end{equation}
Since the corresponding flow hypersurfaces can never touch, in view of the maximum principle, we conclude that the general flow only exists for a finite time and that
\begin{equation}\lae{5.8}
T_1<T^*<T_2,
\end{equation}
where $T_i$ and $T^*$ are the lengths of the corresponding maximal time intervals; for the lower estimate we also used an argument in the proof of \rt{5.6}. 

By the same argument we also obtain:
\bl\lal{5.1}
Let $M(t)$ be a solution of \re{5.1} on a maximal time interval $[0,T^*)$ and represent $M(t)$, for a fixed $t\in [0,T^*)$, as a graph in polar coordinates with center in $x_0\in \hat M(t)$,
\begin{equation}
M(t)=\graph u(t,\cdot)),
\end{equation}
then
\begin{equation}\lae{5.10}
\inf_{M(t)}u\le \Theta(t,T^*)\le \sup_{M(t)}u.
\end{equation}
\el
\bp
The sphere with center $x_0$ and radius $\Theta(t,T^*)$ has to intersect $M(t)$ because of \re{5.8}. Note that, when the relation \re{5.6} is valid at time $t=t_0$, then it is also valid for any $t\ge t_0$ provided the flows exist that long.
\ep
The solution $\Theta=\Theta(t,T^*)$ of \re{5.5} is given by 
\begin{equation}
\Theta =\arccos e^{(t-T^*)},
\end{equation}
since
\begin{equation}
(\log \cos\Theta)'=1.
\end{equation}

Let $\rho_-(t)$ \resp $\rho_+(t)$ be the \tit{inradius} \resp \tit{circumradius} of $\hat M(t)$. Choosing their respective centers as origins of geodesic polar coordinates we deduce from \re{5.10}
\begin{equation}\lae{5.13}
\rho_-(t)\le \Theta (t,T^*)\le \rho_+(t),
\end{equation}
i.e., 
\begin{equation}
\lim_{t\ra T^*}\rho_-(t)=0.
\end{equation}
We want to prove that the corresponding limit of $\rho_+(t)$ also vanishes. Then, the flow would shrink to a point.

Let $x_0\in \hat M(t)$ be arbitrary and consider the corresponding conformally flat coordinate system
\begin{equation}\lae{5.15}
d\bar s^2=\frac1{(1+\frac14 r^2)^2}\{dr^2+r^2\s_{ij}d\xi^id\xi^j\}.
\end{equation}
Write $M(t)$ as graph of $u(t)$ in Euclidean polar coordinates and let $\ka_i$ \resp $\tilde\ka_i$ be the principal curvatures of $M(t)$ when considered as a hypersurface in $\Ss$ \resp $\R$, then we can prove:
\bl
The principal curvatures $\tilde\ka_i$ of $M(t)$ are pinched, i.e., there exists a uniform constant $c$ such that
\begin{equation}\lae{5.16}
\tilde\ka_n\le c \tilde \ka_1,
\end{equation}
where the $\tilde\ka_i$ are labelled
\begin{equation}
\tilde\ka_1\le\cdots\le\tilde\ka_n.
\end{equation}
\el
\bp
The $\ka_i$ and $\tilde\ka_i$ are related through the formula
\begin{equation}\lae{5.18}
\frac1{1+\frac14 r^2}\ka_i=\tilde\ka_i-\tfrac12\frac u{1+\frac14 u^2}v^{-1},
\end{equation}
where
\begin{equation}
v^2=1+u^{-2}\s^{ij}u_iu_j\equiv 1+\abs{Du}^2,
\end{equation}
\cf \cite[equ. (1.1.51)]{cg:cp}. Hence, we deduce
\begin{equation}
\tilde\ka_1\ge \frac1{1+\frac14 c_0^2}\ka_1\ge \frac1{1+\frac14 c_0^2} c_1=c_1',
\end{equation}
since in view of \rl{5.3} below
\begin{equation}\lae{5.21}
u\le c_0,
\end{equation}
where $c_0=c_0(M_0)$ is a uniform constant, and we conclude further
\begin{equation}
\tilde\ka_n\le \ka_n+\tfrac12 c_0
\end{equation}
yielding
\begin{equation}
\frac{\tilde\ka_n}{\tilde\ka_1}\le (1+\frac14 c_0^2)\frac{\ka_n}{\ka_1}+\frac{c_0}{2c_1'}\le (1+\frac14 c_0^2)c_2+\frac{c_0}{2c_1'}
\end{equation}
because of \re{5.4}.
\ep
\bl\lal{5.3}
Let $x_0\in \hat M(t)$ be as above and let $M(t)=\graph u$ be a representation of $M(t)$ in Euclidean polar coordinates, then there exists a constant $c_0=c_0(M_0)$ such that the estimate \re{5.21} is valid for any $t\in [0,T^*)$. Moreover, for any $T\in [0,T^*)$ and $x_0\in \hat M(T)\su \Ss$, the flow hypersurfaces $M(t)$, $0\le t\le T$, can be represented  as graphs in the geodesic polar coordinate system of $\Ss$ with center in $x_0$.
\el
\bp
The convex bodies $\hat M(t)\su \Ss$ are decreasing with respect to  $t$, especially, we have
\begin{equation}\lae{5.24}
\hat M(t)\su \hat M_0\qq \A\, t\in [0,T^*),
\end{equation}
\cf \rr{5.5} below. Since $\hat M_0$ is strictly convex its diameter is less than $\pi$
\begin{equation}\lae{5.25}
\diam M_0<\pi -\ga, \q\ga>0.
\end{equation}
Hence, any geodesic starting in $x_0$ which is contained in $\ol{\hat M(t)}$ has length less than $\pi-\ga$, which in turn implies that the estimate \re{5.21} should be valid with $c_0=c_0(\ga)$.

The second claim of the lemma  is an immediate consequence of \re{5.24} and \re{5.25}.
\ep

Now, choose $x_0\in \hat M(t)$ to be the center of the inball of $\hat M(t)\su \Ss$ with corresponding inradius $\rho_-(t)$ and circumradius $\rho_+(t)$, and let $\tilde\rho_-(t)$ \resp $\tilde\rho_+(t)$ be the inradius \resp circumradius of $\hat M(t)\su \R$. Note that the center of the Euclidean inball is the center of the polar coordinates.

The pinching estimate \re{5.16} then implies, \cf \cite[Theorem 5.1 and Lemma 5.4]{andrews:euclidean},
\begin{equation}\lae{5.26}
\tilde\rho_+(t)\le c\tilde\rho_-(t)
\end{equation}
with a uniform constant $c$, hence $\hat M(t)\su \R$ is contained in the Euclidean ball $B_{\tilde\rho}(0)$
\begin{equation}
\hat M(t)\su B_{\tilde\rho}(0),\qq\tilde\rho=2c\tilde\rho_-(t).
\end{equation}
Define $\tilde \Theta$ by
\begin{equation}
\tilde \Theta=2\tan \frac\Theta2,
\end{equation} 
then we deduce from \re{5.10}
\begin{equation}\lae{5.29}
\inf_{M(t)}u\le \tilde\Theta\le\sup_{M(t)}u,
\end{equation}
where $M(t)=\graph u$ is now a representation of $M(t)$ in Euclidean polar coordinates, concluding further
\begin{equation}
\tilde\rho(t)=2c\tilde\rho_-(t)\le2c \tilde\Theta.
\end{equation}
Choose $\de>0$ so small such that
\begin{equation}\lae{5.31}
2c\tilde\Theta(t,T^*)\le1\qq\A\, \abs{T^*-t}\le\de,
\end{equation}
then
\begin{equation}
\tilde\rho(t)\le1,
\end{equation}
hence, in $\Ss$, we have
\begin{equation}
\hat M(t)\su B_{\rho(t)}(x_0),
\end{equation}
where $B_{\rho(t)}(x_0)$ is the geodesic ball with center $x_0$ and radius
\begin{equation}
\rho(t)=\int_0^{\tilde\rho(t)}\frac1{1+\frac14 r^2}=2\arctan \frac{\tilde\rho(t)}2,
\end{equation}
i.e.,
\begin{equation}
\rho\le\tilde\rho\q\wed\q \rho\ge \frac{\tilde\rho}2.
\end{equation}
Thus, we have proved:
\bl
Let $B_{\rho_-(t)}(x_0)\su \hat M(t)$ be an inball, then
\begin{equation}\lae{5.36}
\hat M(t)\su B_{4c\rho_-(t)}(x_0)\qq\A\, t\in [T^*-\de,T^*),
\end{equation}
where $c$ is the constant in \re{5.26}, or equivalently,
\begin{equation}\lae{5.37}
\rho_+(t)\le 4c\rho_-(t).
\end{equation}
Hence, the flow \re{5.1} converges to a point.
\el
\br\lar{5.5}
The convex bodies $\hat M(t)$ converge monotonely, i.e.,
\begin{equation}
t_1<t_2\q\im\q \hat M(t_2)\su \hat M(t_1),
\end{equation}
yielding
\begin{equation}
p\in \hat M(t)\qq\A\,t\in[0,T^*).
\end{equation}
\er
\bp
It suffices to consider $t_2-t_1$ to be small such that $M(t)$, $t\in [t_1,t_2]$,  can be written as graphs in polar coordinates with center in $\hat M(t_2)$. Then $u=u(t,\cdot)$ satisfies the scalar flow equation
\begin{equation}
\dot u=-Fv<0.
\end{equation}
\ep

Let us finish this section by proving that the flow hypersurfaces are smooth and uniformly convex during the evolution.
\bt\lat{5.6}
During the evolution the flow hypersurfaces $M(t)$ are smooth and uniformly convex satisfying a priori estimates in any compact subinterval
\begin{equation}
[0,T]\su [0,T^*),
\end{equation}
where the a priori estimates only depend on $M_0$, $F$ and $T$.
\et
\bp
It suffices to prove the a priori estimates. Let $0<T<T^*$, then the inradius $\rho_-(t_0)$ satisfies
\begin{equation}\lae{5.42}
0<c\Theta(T,T^*)\le \rho_-(T)
\end{equation}
with a uniform constant independent of  $T$. Indeed, from \re{5.26} and \re{5.29} we infer
\begin{equation}\lae{5.43}
\tilde\Theta(T,T^*)\le c \tilde\rho_-(T),
\end{equation}
where
\begin{equation}
\theta(T,T^*)=\int_0^{\tilde\Theta(T,T^*)}\frac1{1+\frac14r^2}
\end{equation}
and
\begin{equation}
\rho_-(T)=\int_0^{\tilde\rho_-(T)}\frac1{1+\frac14r^2}.
\end{equation}
On the other hand,  $\tilde\rho_-(T)$ as well as $\tilde\Theta(T,T^*)$ are uniformly bounded by the constant $c_0$, in view of \re{5.21} and \re{5.29}. The estimate \re{5.42} is therefore an immediate consequence of \re{5.43}.
Let $x_0\in\hat M(T)$ be the center of an inball and introduce geodesic polar coordinates with center $x_0$. Then, the coordinate system covers the flow (5.1) as long as $0\le t\le T$, in view of \rl{5.3}. Writing the flow hypersurfaces as graphs of a function $u(t,\cdot)$ we have
\begin{equation}\lae{5.46}
0<\de\le u\le \pi-\ga
\end{equation}
and hence,  due to the convexity of  $M(t)$,
\begin{equation}
v^2=1+\sin^{-2}u\,\s^{ij}u_iu_j
\end{equation}
is uniformly bounded. Furthermore, we have already proved that the principal curvatures are uniformly bounded from below
\begin{equation}
0<c_1\le \ka_i.
\end{equation}
Since $F$ is concave it suffices to prove that the $\ka_i$ are also uniformly bounded from above
\begin{equation}\lae{5.49}
\ka_i\le c_2(T)\qq\A\,0\le t\le T
\end{equation}
in order to first apply  the Krylov-Safonov and then the Schauder estimates to obtain the desired a priori estimates.

To derive \re{5.49} we consider the function
\begin{equation}\lae{5.50}
\chi=\frac1{\sin u}v,
\end{equation}
which satisfies the evolution equation \fre{1.5.314}. Let $\tilde\chi=\chi^{-1}$, then $\tilde \chi$ solves the evolution equation
\begin{equation}
\dot{\tilde\chi}-F^{ij}\tilde\chi_{ij}=F^{ij}h_{ki}h^k_j\tilde\chi-2F\frac{\bar H}nv\tilde\chi.
\end{equation}
Because of \re{5.46} and the boundedness of $v$ there exists $\de>0$ such that
\begin{equation}
\tilde\chi>2\de\qq\A\, t\in [0,T]
\end{equation}
and hence
\begin{equation}
\f=\log(\tilde\chi-\de)
\end{equation}
is well defined and satisfies the evolution equation  
\begin{equation}\lae{5.54}
\begin{aligned}
\dot\f-F^{ij}\f_i\f_j=F^{ij}h_{ki}h^k_j\frac{\tilde\chi}{\tilde\chi-\de}+F^{ij}\f_i\f_j-2F\frac{\bar H}nv\frac{\tilde\chi}{\tilde\chi-\de}.
\end{aligned}
\end{equation}
We are now ready to prove the estimate \re{5.49}. As in the proof of \frl{4.1} we may pretend that $h^n_n=\ka_n$, the largest principal curvature, is a smooth function and look at the point $(t_0,\xi_0)$, $t_0>0$, where
\begin{equation}
w=\log h^n_n-\f
\end{equation}
assumes its maximum in $[0,T]\times \Ss[n]$.

Applying the maximum principle we obtain
\begin{equation}
\begin{aligned}
0&\le -F^{ij}h_{ki}h^k_j\frac\de{\tilde\chi-\de} +K_N\{2F(h^n_n)^{-1}-F^{kl}g_{kl}\}\\
&\hp{\le}\;+2F\frac{\bar H}nv\frac{\tilde\chi}{\tilde\chi-\de}.
\end{aligned}
\end{equation}
Since $F^{ij}$ is uniformly positive definite and
\begin{equation}
F\le c h^n_n,
\end{equation}
we deduce $w$ and, hence, $h^n_n$ is a priori bounded.
\ep
\br\lar{5.7}
Let $\de$ be the small constant in \re{5.31} and define
\begin{equation}
t_\de=T^*-\de,
\end{equation}
then we deduce from \re{5.36}
\begin{equation}
\hat M(t_\de)\su B_{8c\rho_-(t_\de)}(x_0)\qq\A\, x_0\in \hat M(t_\de).
\end{equation}
Choosing $\de$ even a bit smaller without changing the notation we may also assume that
\begin{equation}
8c\rho_-(t_\de)\le 8c\Theta(t_\de,T^*)<1.
\end{equation}
In view of the a priori estimates in the preceding theorem we shall henceforth only consider $t\in [t_\de,T^*)$.
\er

\section{The rescaled flow}
We shall first prove that
\begin{equation}
\Theta(t,T^*)\sup_{M(t)}F\le \const\qq\A\,t_\de\le t<T^*.
\end{equation}
The proof will be an adaptation of the proof of a similar result in \cite[Theorem 7.5]{andrews:euclidean}. Let $t_\de<t_0<T^*$ be arbitrary and $B_{\rho_-(t_0)}(x_0)$ be an inball of $\hat M(t_0)$. Choosing $x_0$ to be the center of a geodesic polar coordinate system the hypersurfaces $M(t)$ can be represented as graphs
\begin{equation}
M(t)=\graph u(t,\cdot)\qq\A\,t_\de\le t \le t_0
\end{equation}
such that 
\begin{equation}\lae{6.3}
\rho_-(t_0)\le u(t_0)\le u(t)\le 1,
\end{equation}
\cf \rr{5.7}.
\bl\lal{6.1}
Let $\chi$ be defined as in \fre{5.50}, then
\begin{equation}
\chi_i=0\q\im\q u_i=0.
\end{equation}
\el
\bp
The function
\begin{equation}
\h(r)=\frac1{\sin r}
\end{equation}
is a solution of the equation
\begin{equation}
\dot\h=-\frac{\bar H}n\h,
\end{equation}
where $\bar H$ is the mean curvature of the slices $\{x^0=r\}$. Moreover,
\begin{equation}
v^{-2}=1-\norm{Du}^2
\end{equation}
implies
\begin{equation}
v_i=-h_{ij}u^jv^2+\frac{\bar H}nu_iv,
\end{equation}
where
\begin{equation}
u^j=g^{ij}u_i\q\wed\q \norm{Du}^2=g^{ij}u_iu_j.
\end{equation}
On the other hand, we deduce from
\begin{equation}
0=\chi_i=\dot\h u_iv+\h v_i=-\frac{\bar H}n\h u_i v+\h v_i
\end{equation}
\begin{equation}
v_i=\frac{\bar H}nu_iv
\end{equation}
concluding further
\begin{equation}
h_{ij}u^j=0
\end{equation}
and thus $u_i=0$, since $h_{ij}$ is positive definite.
\ep

Let $\tilde \chi=\chi^{-1}$ as before, then $\tilde\chi$ is the equivalent of the Euclidean support function and in view of the estimate \re{6.3} and \rl{6.1} there exists a universal constant $\e_0$ such that
\begin{equation}\lae{6.13}
0<\tilde\chi-2\e_0\rho_-(t_0)\qq\A\,t_\de\le t\le t_0.
\end{equation}
We are now able to prove:
\bl
There exists a uniform constant $c$ such that
\begin{equation}\lae{6.14}
\Theta(t,T^*)F\le c\qq\A\, t_\de\le t <T^*.
\end{equation}
\el
\bp
Let $t_0\in (t_\de, T^*)$ be arbitrary and consider the function
\begin{equation}
\f=\log F-\log (\tilde \chi-\e_0\rho_-(t_0))
\end{equation}
in the interval $[t_\de, t_0]$. Define
\begin{equation}
w(t)=\sup_{M(t)}\f,
\end{equation}
then $w$ satisfies the differential inequality, \cf \fre{5.54},
\begin{equation}
\begin{aligned}
\dot w&\le -F^{ij}h_{ki}h^k_j \frac{\e_0\rho_-(t_0)}{\tilde\chi-\e_0\rho_-(t_0)}+K_NF^{ij}g_{ij}\\
&\hp{\le}\;+2F\frac{\bar H}nv\frac{\tilde\chi}{\tilde\chi-\e_0\rho_-(t_0)}\\
&\le -\frac1{F^{ij}g_{ij}}F^2\frac{\e_0\rho_-(t_0)}{\tilde\chi-\e_0\rho_-(t_0)}+cK_N+cF\frac1{\tilde\chi-\e_0\rho_-(t_0)},
\end{aligned}
\end{equation}
where we used that
\begin{equation}
\frac{\bar H}n=\frac{\cos u}{\sin u}
\end{equation}
and 
\begin{equation}\lae{6.19}
\tilde\chi=\sin u v^{-1}.
\end{equation}
Setting 
\begin{equation}
\tilde w=e^w
\end{equation}
we infer
\begin{equation}
\begin{aligned}
\dot{\tilde w}&\le -\tilde c \tilde w^3\e_0\rho_-(t_0)(\tilde\chi-\e_0\rho_-(t_0))+cK_N\tilde w+c \tilde w^2\\
&\le \tilde w^2\{c+cK_N\tilde w^{-1}-\tilde c \e_0^2\rho_-(t_0)^2\tilde w\},
\end{aligned}
\end{equation}
in view of \re{6.13}.

Hence we conclude
\begin{equation}
\sup_t\tilde w\le \max(\tilde w(t_\de)+1, c\e_0^{-2}\rho_-(t_0)^{-2}),
\end{equation}
where $c$ is a new uniform constant independent of $t_0$ and $t_\de$. Choosing $t_0$ large enough we obtain
\begin{equation}
\tilde w(t_0)\rho^2_-(t_0)\le c \e_0^{-2}
\end{equation}
and thus, because of \re{6.13}, 
\begin{equation}
\rho_-(t_0)\sup_{M(t_0)}F\le c \e_0^{-2}
\end{equation}
with a different constant $c$. To complete the proof we use the estimates \fre{5.42}.
\ep
\bc
The rescaled principal curvatures $\tilde\ka_i=\Theta\ka_i$ satisfy
\begin{equation}\lae{6.25}
\tilde\ka_i\le c\qq\A\,t_\de\le t<T^*
\end{equation}
with a uniform constant.
\ec
\bp
From \re{6.14} we infer
\begin{equation}
c\ge \tilde F=F\Theta=\sum_iF^i_i\tilde\ka_i.
\end{equation}
Since $F^i_j$ is uniformly positive definite because of the pinching estimates, the result follows.
\ep

Next we want to apply the Harnack inequality to get an estimate from below for $\tilde F$
\begin{equation}
\inf_{M(t)}\tilde F\ge c>0\qq\A\, t_\de\le t<T^*.
\end{equation}
To convince ourselves that the necessary requirements are fulfilled we first have to establish some preparatory results.
\bl\lal{6.4}
Let $t_1\in [t_\de,T^*)$ be arbitrary and let $t_2>t_1$ be such that
\begin{equation}\lae{6.28}
\Theta(t_2,T^*)=\tfrac12\Theta(t_1,T^*).
\end{equation}
Let $x_0\in\hat M(t_2)$ be the center of an inball. Introduce polar coordinates around $x_0$ and write the hypersurfaces $M(t)$ as graphs
\begin{equation}
M(t)=\graph u(t,\cdot),
\end{equation}
then there exists a positive constant c such that
\begin{equation}\lae{6.30}
c^{-1}\Theta(t_2,T^*)\le u(t,\xi)\le c\Theta(t_2,T^*)\qq\A\, t\in [t_1,t_2],
\end{equation}
and hence
\begin{equation}\lae{6.31}
\frac{u_{\tup{max}}(t)}{u_{\tup{min}}(t)}\le c^2\qq\A\, t\in [t_1,t_2],
\end{equation}
where
\begin{equation}
u_{\tup{max}}(t)=\sup_{M(t)}u\q\wed\q u_{\tup{min}}(t)=\inf_{M(t)}u.
\end{equation}
\el
\bp
Let $B_{\rho_-(t_1)}(y_0)$ be an inball of $\hat M(t_1)$, then we infer from \fre{5.36} and \fre{5.13}
\begin{equation}
\begin{aligned}\lae{6.33}
\hat M(t_1)\su B_{4c\rho_-(t_1)}(y_0)\su B_{4c\Theta(t_1,T^*)}(y_0)\su B_{8c\Theta(t_2,T^*)}(y_0)
\end{aligned}
\end{equation}
and we deduce further, since
\begin{equation}
\hat M(t_2)\su \hat M(t_1),
\end{equation}
\begin{equation}
\hat M(t_1)\su B_{16c\Theta(t_2,T^*)}(x_0).
\end{equation}
Hence, we have proved the upper estimate in \re{6.30}. The lower estimate follows from \re{5.37} and \re{5.13}, because
\begin{equation}
\rho_-(t_2)\le u(t,\xi)\qq\A\, t\in [t_1,t_2].
\end{equation}
\ep

\bl\lal{6.5}
Under the assumptions of the preceding lemma the quantity
\begin{equation}
v^2=1+\sin^{-2}u\,\s^{ij}u_iu_j
\end{equation}
is uniformly bounded in $[t_1,t_2]\times \Ss[n]$.
\el
\bp
From \cite[inequality (2.7.83)]{cg:cp} we obtain
\begin{equation}
v(t,\xi)\le e^{\bar\ka(u_{\tup{max}}-u_{\tup{min}})},
\end{equation}
where $0\le\bar\ka$ is an upper bound for the principle curvatures of the slices $\{x^0=\const\}$ intersecting $M(t)$, hence
\begin{equation}
\bar\ka\le \frac1{\sin{u_{\tup{min}}}}\le c\frac1{u_{\tup{min}}}
\end{equation}
and
\begin{equation}
v(t,\xi)\le e^{c\big(\frac{u_{\tup{max}}}{u_{\tup{min}}}-1\big)}.
\end{equation}
Combining this estimate with the one in \re{6.31} gives the result.
\ep

\bl\lal{6.6}
Define $\vt$ by
\begin{equation}
\vt(r)=\sin r
\end{equation}
and
\begin{equation}\lae{6.42}
\begin{aligned}
\f&=\int_{r_2}^u\vt^{-1}\\
&=\{\log(\sin\tfrac r2)-\log(\cos \tfrac r2)\}_{\big|^u_{r_2}},
\end{aligned}
\end{equation}
where, $r_2=\Theta(t_2,T^*)$, then $\f(t,\cdot)$ is uniformly bounded in $C^2(\Ss[n])$ for any $t_1\le t\le t_2$, independent of $t_1$, $t_2$. Furthermore, let $\ch ijk$ \resp $\tilde{\ch ijk}$ be the Christoffel symbols of the metrics $g_{ij}$ \resp $\s_{ij}$, then the tensor
\begin{equation}\lae{6.43}
\ch ijk-\tilde{\ch ijk}
\end{equation}
is also uniformly bounded independent of $t_1$, $t_2$.
\el
\bp
The $C^0$ and the $C^1$-estimates are due to \re{6.30} and \rl{6.5}. To prove the $C^2$-estimates we employ the relation
\begin{equation}\lae{6.44}
h^i_j=v^{-1}\vt^{-1}\{-(\s^{ik}-v^{-2}\f^i\f^k)\f_{jk}+\dot\vt\de^i_j\},
\end{equation}
\cf \cite[equ. (3.26)]{cg:icf-hyperbolic}, where
\begin{equation}
\f^i=\s^{ik}\f_k
\end{equation}
and where the covariant derivatives are with respect to the metric $\s_{ij}$. Multiplying both sides of \re{6.44} with $\Theta(t,T^*)$ we deduce
\begin{equation}
\norm{\f_{ij}}\le c\qq\A\, t\in [t_1,t_2],
\end{equation}
in view of the $C^1$-estimates, \re{6.30} and \re{6.25}.

To prove the boundedness of \re{6.43} we choose coordinates such that in a fixed point $\tilde{\ch ijk}$ vanishes. Then $\ch ijk$ is a uniformly bounded tensor comprised of algebraic compositions of $v$, $D\f$, $D^2\f$ and $\s_{ij}$ as one easily checks.
\ep

Let us define a new time parameter
\begin{equation}\lae{6.47}
\tau=-\log\Theta,
\end{equation}
then
\begin{equation}
\df t\tau=-\frac\Theta{\dot\Theta}=\Theta \frac{\sin\Theta}{\cos\Theta}.
\end{equation}
Let a prime indicate differentiation with respect to $\tau$ and a dot with respect to $t$, and let us denote scaled quantities by a tilde unless otherwise specified, e.g., let
\begin{equation}
\tilde F=F\Theta,
\end{equation}
then
\begin{equation}\lae{6.50}
\tilde F'=\dot F \Theta^2 \frac{\sin\Theta}{\cos\Theta}-\tilde F
\end{equation}
and we shall prove:
\bl
$\tilde F$ satisfies a uniformly parabolic equation of the form
\begin{equation}
\tilde F'-a^{ij}\tilde F_{ij}+b^i\tilde F_i+c\tilde F=0
\end{equation}
in the cylinder
\begin{equation}
Q(\tau_1,\tau_2)=[\tau_1,\tau_2]\times\Ss[n],
\end{equation}
where
\begin{equation}
\tau_i=-\log\Theta(t_i,T^*),
\end{equation}
with uniformly bounded coefficients, and where the covariant derivatives are with respect to standard metric $\s_{ij}$ of $\Ss[n]$. The coefficients are bounded independently of $\tau_i$. Since, in view of \re{6.28}
\begin{equation}
\tau_2=\tau_1+\log2,
\end{equation}
we deduce, by applying the parabolic Harnack inequality,
\begin{equation}
\sup_{M(t_1)}\tilde F\le c\inf_{M(t_2)}\tilde F
\end{equation}
with a uniform constant $c$.
\el
\bp
It suffices to prove that $\tilde F$ satisfies a uniformly parabolic equations as indicated. Combining \re{6.50} and \fre{1.2.314} we immediately deduce, in view of \re{6.25} and the pinching estimates, that the only non-trivial term in the transformation of \re{1.2.314} is
\begin{equation}
-F^{ij}F_{;ij}\Theta^2\frac{\sin\Theta}{\cos\Theta},
\end{equation}
where the semicolon indicates covariant derivatives with respect to $g_{ij}$.

Now, using geodesic polar coordinates as in \rl{6.4}, we can express the metric in the form
\begin{equation}
g_{ij}=\sin^2u(\f_i\f_j+\s_{ij}),
\end{equation}
\cf the definition in \re{6.42}, and we deduce
\begin{equation}
g^{ij}\Theta^2
\end{equation}
is uniformly positive definite, in view of \re{6.30} and \rl{6.5}, hence
\begin{equation}
\Theta^2 F^{ij}=F^i_kg^{kj}\Theta^2
\end{equation}
is uniformly positive definite.

Thus, it remains to consider the covariant derivatives, but
\begin{equation}
F_{;ij}=F_{ij}-\{\ch ijk-\tilde{\ch ijk}\}F_k,
\end{equation}
where $F_{ij}$ are the covariant derivatives of $F$ with respect to $\s_{ij}$ and $\ch ijk$ \resp $\tilde{\ch ijk}$ are the Christoffel symbols with respect to $g_{ij}$ \resp $\s_{ij}$, hence we infer from \rl{6.6}
\begin{equation}
-F^{ij}F_{;ij}\Theta^2\frac{\sin\Theta}{\cos\Theta}=-a^{ij}\tilde F_{ij}+b^i\tilde F_i,
\end{equation}
where $a^{ij}$ is uniformly positive definite and $b^i$ uniformly bounded.
\ep

\bc\lac{6.8}
The scaled curvatures $\tilde\ka_i$ are uniformly bounded from below
\begin{equation}
\tilde\ka_i=\ka_i\Theta\ge c>0.
\end{equation}
\ec
\bp
Since
\begin{equation}
\inf_{M(t)}\tilde F\le \inf_{M(t)}\tilde\ka_n\le c\inf_{M(t)}\tilde\ka_1,
\end{equation}
where the $\tilde\ka_i$ are labelled
\begin{equation}
\tilde\ka_1\le\cdots\le \tilde\ka_n,
\end{equation}
and $t_1\in[t_\de,T^*)$ is arbitrary, it suffices to estimate
\begin{equation}
\sup_{M(t)}\tilde F\ge c>0\qq\A\,t_1\le t\le t_2.
\end{equation}

Indeed, let $(t,\xi)\in M(t)$ be a point such that
\begin{equation}
u(t,\xi)=\sup_{M(t)}u,
\end{equation}
then
\begin{equation}
\ka_i\ge\frac{\cos u}{\sin u}
\end{equation}
and
\begin{equation}
\tilde\ka_i\ge \frac{\cos u}{\sin u}\Theta\ge c>0,
\end{equation}
because of \re{6.30} and hence
\begin{equation}
\sup_{M(t)}\tilde F\ge F(\tilde\ka_i(t,\xi))\ge c>0.
\end{equation}
\ep

Now, let $x_0\in\Ss$ be the point the flow hypersurfaces are shrinking to and introduce geodesic polar coordinates around it. Let
\begin{equation}
M(t)=\graph u(t,\cdot)
\end{equation}
and let
\begin{equation}
\tilde u(\tau,\xi)=u(t,\xi)\Theta(t,T^*)^{-1},
\end{equation}
where $\tau$ is defined as in \re{6.47}. Then, we can prove:
\bl
There exists a uniform constant $c$ such that
\begin{equation}
\tilde u\ge c>0\qq\A\,\tau \in  Q(\tau_\de,\un),
\end{equation}
where
\begin{equation}\lae{6.73}
\tau_\de=-\log(\Theta(t_\de,T^*))
\end{equation}
and
\begin{equation}
Q(\tau_\de,\un)=[\tau_\de,\un)\times \Ss[n].
\end{equation}
\el
\bp
Let us look at the rescaled version of the scalar curvature equation
\begin{equation}
\dot u=\pde ut=-Fv,
\end{equation}
which has the form
\begin{equation}\lae{6.76} 
\begin{aligned}
\tilde u'&=\dot u\frac{\sin\Theta}{\cos\Theta}+\tilde u\\
&=-\tilde F \Theta^{-1}\frac{\sin\Theta}{\cos\Theta}v+\tilde u\\
&\le -2c+\tilde u,
\end{aligned}
\end{equation}
in view of \rc{6.8}.

Let us suppose there exists $\tau_0\ge\tau_\de$ and $\xi\in \Ss[n]$ such that
\begin{equation}
\tilde u(\tau_0,\xi)\le c,
\end{equation}
then
\begin{equation}
\tilde u'\le -c\qq\A\, \tau\ge \tau_0,
\end{equation}
where $\tilde u$ is evaluated at $(\tau,\xi)$, yielding
\begin{equation}
\tilde u(\tau)-\tilde u(\tau_0)\le -c(\tau-\tau_0)\qq\A\, \tau_0\le\tau<\un,
\end{equation}
a contradiction, hence we conclude
\begin{equation}
\tilde u\ge c\qq\A\, (\tau,\xi)\in  Q(\tau_\de,\un).
\end{equation}
\ep

\bl
The quantities $\tilde u$, $v$ and $\abs{D\tilde u}$ are uniformly bounded in $Q(\tau_\de,\un)$, where
\begin{equation}
\abs{D\tilde u}^2=\s^{ij}\tilde u_i\tilde u_j.
\end{equation}
\el
\bp
(i) Let $t\in [t_\de,T^*)$ be arbitrary, and $B_{\rho_-(t)}(y_0)$ be an inball of $\hat M(t)$, then we infer from \re{6.33}
\begin{equation}
\hat M(t)\su B_{4c\rho_-(t)}(y_0).
\end{equation}
On the other hand, $x_0\in\hat M(t)$ and
\begin{equation}
\rho_-(t)\le \Theta(t,T^*),
\end{equation}
hence
\begin{equation}
\hat M(t)\su B_{8c\Theta(t,T^*)}(x_0)
\end{equation}
yielding
\begin{equation}
\tilde u\le 8c.
\end{equation}

(ii) From the proof of \rl{6.5} we immediately deduce that
\begin{equation}
v^2=1+\frac1{\sin^2u}\s^{ij}u_iu_j\le c
\end{equation}
which in turn implies
\begin{equation}
\s^{ij}\tilde u_i\tilde u_j\le c \sin^2u\Theta^{-2}\le \const.
\end{equation}
\ep
\br
Let $\f$ be such that
\begin{equation}
\f_i=\frac1{\sin u}u_i,
\end{equation}
then the covariant derivatives of $\tilde u$ \resp $\f$ with respect to $\s_{ij}$ satisfy the pointwise estimate
\begin{equation}
\norm{\tilde u_{ij}}\le c\norm{\f_{ij}},
\end{equation}
hence we conclude that the $C^2$-norm of $\tilde u$ is uniformly bounded and also the difference of the Christoffel symbols
\begin{equation}
\ch ijk-\tilde{\ch ijk},
\end{equation}
\cf \rl{6.6} and its proof. Moreover, observing that
\begin{equation}\lae{7.91.1}
\frac{\sin u}u=\vt(u)=\vt(\tilde u e^{-\tau})\ge c_0>0,
\end{equation}
where $\vt$ is a smooth function such that
\begin{equation}\lae{7.92.1}
c_0\le \vt\le c_0^{-1}\qq\A\, t\in [t_\de,T^*),
\end{equation}
and taking a similar estimate for $\cos u$ into account, we conclude from \re{6.44}
\begin{equation}
\frac{\sin\Theta}{\cos\Theta}Fv=F(\tfrac{\sin\Theta}{\cos\Theta}h^i_j)v=\F(x,e^{-\tau},\tilde u,\tilde ue^{-\tau}, D\tilde u,D^2\tilde u),
\end{equation}
where $\F$ is smooth function with respect to its arguments, monotone and concave with respect to $-\tilde u_{ij}$, where the covariant derivatives are defined relative to the standard metric on $\Ss[n]$.
\er
Hence, we deduce, by applying the Krylov-Safonov and Schauder estimates:
\bt\lat{6.12}
The rescaled function $\tilde u$ satisfies the uniformly parabolic equation
\begin{equation}
\tilde u'=-\F+\tilde u
\end{equation}
in $Q(\tau_\de,T^*)$ and $\tilde u(\tau,\cdot)$ obeys uniform a priori estimates in $C^\un(\Ss[n])$ independently of $\tau$.
\et

In the next  section we shall prove that $\tilde u$ converges exponentially fast to the constant function $1$ when $F$ is strictly concave or when $F=\tfrac1n H$.

Let us also emphasize that the $\Theta\ka_i$ are not the principal curvatures of $\graph\tilde u$, though they are of course related.

\section{Convergence to a sphere}

The key estimate for proving that the rescaled hypersurfaces converge to a sphere is the exponential decay of the quantity
\begin{equation}
\abs{\tilde A}^2-\tfrac1n \tilde H^2=\tfrac1n\sum_{i<j}\abs{\tilde\ka_i-\tilde\ka_j}^2.
\end{equation}
Huisken proved it in \cite[Section 5]{huisken:euclidean} by deriving the uniform estimate
\begin{equation}\lae{7.2}
H^{-(2-\s)}\{\abs{A}^2-\tfrac1n H^2\}\le \const\qq\A\, t\in [0,T^*),
\end{equation}
for the unscaled hypersurfaces, where $0<\s<1$ is small.

We shall adapt his approach to the present situation where the fact that we consider general curvature function $F$ creates some additional difficulties. Some of the estimates, we shall prove below, will be valid for arbitrary curvature functions, or at least for curvature functions we consider in this paper, but the estimate \re{7.2} can only be proved for $F=\tfrac1n H$ or $F$  strictly concave. 
\bl
Let $M$ be a strictly convex hypersurface with pinched principal curvatures such that
\begin{equation}
h_{ij}\ge \e_0Hg_{ij},\qq \e_0>0,
\end{equation}
and let $F$ be monotone and concave. Then there exists $\e>0$, $\e=\e(\e_0,F)$, such that
\begin{equation}\lae{7.9}
\begin{aligned}
Z=Fh^k_ih_{kj}h^{ij}-\abs{A}^2F^{ij}h_{ki}h^k_j\ge 2\e^2H^2\sum_{i<j}\abs{\ka_i-\ka_j}^2,
\end{aligned}
\end{equation}
or equivalently,
\begin{equation}\lae{7.10}
Z\ge 2\e^2 F^{ij}h_{ki}h^k_j\sum_{i<j}(\ka_i-\ka_j)^2.
\end{equation}
\el
\bp
Huisken proved the lemma for $F=H$. We consider $F$ to be defined in $\C_+$ and set
\begin{equation}
F_i=\pd F{\ka}i.
\end{equation}
Let us also label the $\ka_i$ such that
\begin{equation}
\ka_1\le\cdots\le \ka_n,
\end{equation}
then
\begin{equation}
F_1\ge \cdots\ge F_n,
\end{equation}
because $F$ is concave. Writing
\begin{equation}
\sum_{i\not=j}F_i\ka_i\ka_j^3=\sum_{i<j}F_i\ka_i\ka_j^3+\sum_{j<i}F_i\ka_i\ka_j^3
\end{equation}
and
\begin{equation}
-\sum_{i\not=j}F_i\ka_i^2\ka_j^2=-\sum_{i<j}F_i\ka_i^2\ka_j^2-\sum_{j<i}F_i\ka_i^2\ka_j^2
\end{equation}
we deduce from \re{7.9}
\begin{equation}
\begin{aligned}
Z&=\sum_{i<j}F_i\ka_i\ka_j(\ka_j^2-\ka_i\ka_j)+\sum_{j<i}F_i\ka_i\ka_j(\ka_j^2-\ka_i\ka_j)\\
&=\sum_{i<j}F_i\ka_i\ka_j(\ka_j^2-\ka_i\ka_j)+\sum_{i<j}F_j\ka_i\ka_j(\ka_i^2-\ka_i\ka_j)\\
&\ge \sum_{i<j}F_j\ka_i\ka_j(\ka_i-\ka_j)^2\ge \e H^2\sum_{i<j}(\ka_i-\ka_j)^2.
\end{aligned}
\end{equation}
\ep
Since $F$ is concave satisfying $F(1,\ldots,1)=1$ we have
\begin{equation}\lae{7.17.1}
F\le \tfrac1nH,
\end{equation}
hence
\begin{equation}\lae{7.17}
\abs{A}^2-\tfrac1nH^2\le \abs{A}^2-nF^2.
\end{equation}
We also need a reverse inequality:
\bl\lal{7.4}
Under the assumptions of the previous lemma there exists a positive constant $c$ such that
\begin{equation}
\abs{A}^2-nF^2\le c(\abs{A}^2-\tfrac1nH^2).
\end{equation}
\el
\bp
The proof will reveal that curvatures need not be positive, it will only be necessary that
\begin{equation}
\frac{\ka_i}{\abs A}
\end{equation}
are compactly contained in the defining cone. To simplify the notation we shall also assume that $F(1,\ldots,1)=n$ such that we have to prove the inequality
\begin{equation}
\abs A^2-\tfrac1nF^2\le c(\abs A^2-\tfrac1nH^2),
\end{equation}
or equivalently,
\begin{equation}\lae{7.21}
H^2-F^2\le c\sum_{i<j}(\ka_i-\ka_j)^2.
\end{equation}
Let
\begin{equation}
\f=H^2-F^2
\end{equation}
and consider the convex combination
\begin{equation}
\ka_i(t)=(1-t)\ka_n+t\ka_i,
\end{equation}
where the $\ka_i$ are labelled such that
\begin{equation}
\ka_1\le\cdots\le \ka_n.
\end{equation}
Denote the partial derivatives of $\f$ simply by indices, then
\begin{equation}
\f(\ka_n,\dots,\ka_n)=0\q\wed\q \f_i(\ka_n,\dots,\ka_n)=0,
\end{equation}
hence we deduce from Taylor's formula
\begin{equation}
\f(\ka_i)=\tfrac12\f_{ij}(\ka_i(t))(\ka^i-\ka^n)(\ka^j-\ka^n)
\end{equation}
for some $0\le t\le 1$ yielding the estimate \re{7.21}, since
\begin{equation}
\f_{ij}=2H_iH_j-2F_iF_j-2FF_{ij}
\end{equation}
is uniformly bounded.
\ep

We are going to estimate the function
\begin{equation}\lae{7.28}
f_\s=F^{-\al}(\abs A^2-n F^2),
\end{equation}
where
\begin{equation}
\al=2-\s
\end{equation}
and $0<\s<1$. We shall also drop the subscript $\s$ simply writing $f$ for the left-hand side of \re{7.28}.

In order to derive the evolution equation for $f$ we use the relation
\begin{equation}
f=\abs A^2F^{-\al}-nF^{2-\al}
\end{equation}
and the equations
\begin{equation}
\begin{aligned}
(\abs A^2)'-F^{ij}\abs A^2_{ij}&=2 F^{ij}h_{ki}h^k_j\abs A^2-2 F^{ij}h_{kl;i}h^{kl}_{\hp{kl};j}\\
&\hp{=}\;+2F^{kl.rs}h_{kl;i}h_{rs;j}h^{ij}\\
&\hp{=}\;+4K_NFH-2K_NF^{kl}g_{kl}\abs A^2,
\end{aligned}
\end{equation}
\begin{equation}
\begin{aligned}
(F^{-\al})'-F^{ij}F^{-\al}_{ij}&=-\al F^{ij}h_{ki}h^k_jF^{-\al}-\al(\al+1)F^{-\al-2}F^{ij}F_iF_j\\
&\hp{=}\;-\al K_NF^{kl}g_{kl}F^{-\al} 
\end{aligned}
\end{equation}
and
\begin{equation}
\begin{aligned}
(F^{2-\al})'-F^{ij}F^{2-\al}_{ij}&=-(\al-2)F^{ij}h_{ki}h^k_jF^{2-\al}\\
&\hp{=}\;-(\al-2)(\al-1)F^{ij}F_iF_jF^{-\al}\\
&\hp{=}\;+(2-\al)K_NF^{kl}g_{kl}F^{2-\al}.
\end{aligned}
\end{equation}
We then obtain
\begin{equation}\lae{7.34}
\begin{aligned}
\msp[100]f'-F^{ij}f_{ij}&=\\
&\msp[-200]\s F^{ij}h_{ki}h^k_jf-2F^{ij}\{h_{kl;i}F-h_{kl}F_i\}\{h^{kl}_{\hp{kl};j}F-h^{kl}F_j\}F^{-(2+\al)}\\
&\msp[-204]-\s(1-\s)F^{ij}F_iF_jF^{-2}f+2(\al-1)F^{-1}F^{ij}F_if_j\\
&\msp[-204]+4K_N\{HF-F^{kl}g_{kl}\abs A^2\}F^{-\al}+\s K_NF^{kl}g_{kl}f\\
&\msp[-204]+2F^{kl,rs}h_{kl;i}h_{rs;j}h^{ij}F^{-\al},
\end{aligned}
\end{equation}
where we used the relation
\begin{equation}
\begin{aligned}
F^{ij}\{h_{kl;i}F-h_{kl}F_i\}\{h^{kl}_{\hp{kl};j}F-h^{kl}F_j\}&=F^{ij}h_{kl;i}h^{kl}_{\hp{kl};j}F^2\\
&\msp[-150]-F^{ij}\abs A^2_iF_jF+F^{ij}F_iF_j\abs A^2.
\end{aligned}
\end{equation}

We also need a purely elliptic version of equation \re{7.34}. This can be achieved by replacing $f'$ using the formula
\begin{equation}
\dot h^i_j=F^j_{;i}+Fh^k_ih^j_k+K_NF\de^j_i,
\end{equation}
\cf \cite[Lemma 2.3.3]{cg:cp}. Hence we deduce
\begin{equation}
\begin{aligned}
f'&=2F^{-\al}Z+\s F^{ij}h_{ki}h^k_jf+(2-\al)F^{ij}F_{;ij}F^{-1}f\\
&-2F^{ij}F_{;ij}\abs A^2F^{-(1+\al)}-\al F^{-1}F^{ij}F_{;ij}f\\
&+2\{h^{ij}F^{-\al}-F^{1-\al}nF^{ij}\}F_{;ij}+2h^{ij}F_{;ij}F^{-\al}\\
&+2K_NFHF^{-\al}-\al K_NF^{ij}g_{ij}\abs A^2F^{-\al}\\
&-(2-\al)K_NF^{ij}g_{ij}nF^{2-\al}
\end{aligned}
\end{equation}
concluding further
\begin{equation}\lae{7.38}
\begin{aligned}
\msp[40]-F^{ij}f_{ij}+2F^{-\al}Z&=\\
&\msp[-200]\al F^{ij}F_{;ij}F^{-1}f-2\{h^{ij}-FnF^{ij}\}F_{;ij}F^{-\al}-\s(1-\s)F^{ij}F_iF_jF^{-2}f\\
&\msp[-204]-2F^{ij}\{h_{kl;i}F-h_{kl}F_i\}\{h^{kl}_{\hp{kl};j}F-h^{kl}F_j\}F^{-(2+\al)}\\
&\msp[-204]+2(\al-1)F^{-1}F^{ij}F_if_j+2K_N\{FH-F^{kl}g_{kl}\abs A^2\}F^{-\al}\\
&\msp[-204]+2F^{kl.rs}h_{kl;i}h_{rs;j}h^{ij}F^{-\al}.
\end{aligned}
\end{equation}
Some of the negative terms on the right-hand side can be exploited. First, we observe that
\begin{equation}\lae{7.40}
FH-F^{kl}g_{kl}\abs A^2\le \tfrac1nH^2-\abs A^2=-\tfrac1n \sum_{i<j}(\ka_i-\ka_j)^2.
\end{equation}
In case $F=\tfrac1nH$ it is proved in \cite[Lemma 2.3 (ii)]{huisken:euclidean} that
\begin{equation}\lae{7.41.1}
F^{ij}\{h_{kl;i}F-h_{kl}F_i\}\{h^{kl}_{\hp{kl};j}F-h^{kl}F_j\}\ge \frac1{2n^3}\e^2H^2\abs{DH}^2.
\end{equation}
For more general curvature functions this inequality can not be derived. Instead we shall consider the last term on the right-hand side of \re{7.38}. If $F=F(\ka)$ is strictly concave in a convex cone $\C\su \R[n]$, then there exists a positive constant $c$ such that
\begin{equation}\lae{7.41}
F^{kl,rs}h_{kl;i}h_{rs;j}\le -c \abs A^{-1}\abs{DA}^2g_{ij}
\end{equation}
provided the normalized vectors
\begin{equation}\lae{7.42}
\abs A^{-1}\ka
\end{equation}
stay in a compact set $K\su\C$. The constant then depends on $F$ and $K$. The estimate \re{7.41} was proved in \cite[Lemma 7.12]{andrews:euclidean}.

In our case the principal curvatures of the flow hypersurfaces are pinched, hence, the normalized curvatures \re{7.42} are compactly contained in $\C_+$, and we can prove:
\bl\lal{7.5}
Let the curvature function $F$ satisfy our general assumptions and assume in addition that it is strictly concave, then there exists a uniform constant $\e>0$ such that
\begin{equation}\lae{7.44}
\begin{aligned}
-F^{ij}f_{ij}+2\e^2 F^{ij}h_{ki}h^k_jf&\le \al F^{-1}F^{ij}F_{;ij}f+2(\al-1)F^{-1}F^{ij}F_if_j\\
&\msp[-100]-2\{h^{ij}-FnF^{ij}\}F^{-\al}F_{;ij} -2\e^2\abs{DA}^2F^{-\al}.
\end{aligned}
\end{equation}
\el
\bp
The claim immediately follows from \re{7.10}, \re{7.17}, \re{7.21}, \re{7.40} and \re{7.41} and the fact that $F$ is strictly concave.
\ep
\bl
There exists a uniform constant $c>0$  such that
\begin{equation}\lae{7.49}
\norm{h^{ij}-FnF^{ij}}^2\le c \sum_{i<j}(\ka_i-\ka_j)^2.
\end{equation}
\el
\bp
We have
\begin{equation}
\begin{aligned}
h^{ij}-FnF^{ij}&=\{h^{ij}-\tfrac1nHg^{ij}\}+\{\tfrac1nH-F\}g^{ij}\\
&\hp{=}\;+F(g^{ij}-nF^{ij})\\
&\equiv I_1+I_2+I_3,
\end{aligned}
\end{equation}
where each term can be estimated by the square root of right-hand side of \re{7.49}.

The estimate for $I_1$ is trivial, $I_2$ can be estimated along the lines of the proof of \rl{7.4}, while
\begin{equation}
I_3=Fn(F^{ij}(\ka_n,\ldots,\ka_n)-F^{ij}(\ka_i))
\end{equation}
from which the estimate follows immediately.
\ep
We are now able to prove a crucial estimate:
\bl\lal{7.8}
Let $F$ be strictly concave, then there exists a constant $c>0$ such that for any $p\ge 2$, any $\de>0$ and any $0\le t<T^*$ the estimate
\begin{equation}\lae{7.52}
\begin{aligned}
\e^2\int_MF^{ij}h_{ki}h^k_jf^p&\le \{\de^{-1}c(p-1)+c\}\int_MF^{ij}f_if_jf^{p-2}\\
&\hp{\le}\;+\{\de c(p-1)+c\}\int_M\abs{DA}^2F^{-\al}f^{p-1}
\end{aligned}
\end{equation}
is valid.
\el
\bp
Multiplying inequality \re{7.44} with $f^{p-1}$ and integrating by parts we obtain
\begin{equation}
\begin{aligned}
(p-1)\int_MF^{ij}f_if_jf^{p-2}+2\e^2\int_MF^{ij}h_{ki}h^k_jf^p&\le\\
&\msp[-300]\int_MF^{ij,kl}h_{kl;j}f_if^{p-1}+\al\int_MF^{-1}F^{ij}F_{;ij}f^p\\
&\msp[-400]-2\int_M\{h^{ij}-FnF^{ij}\}F^{-\al}F_{;ij}f^{p-1}+2(\al-1)\int_MF^{-1}F^{ij}F_if_jf^{p-1}\\
&\msp[-300]-2\e^2\int_M\abs{DA}^2f^{-\al}f^{p-1}.
\end{aligned}
\end{equation}
The terms on the right-hand side can be estimated or transformed as follows:
\begin{equation}
\begin{aligned}
\int_MF^{ij,kl}h_{kl;j}f_if^{p-1}&\le \de^{-1}(p-1)\int_MF^{ij}f_if_jf^{p-2}\\
&\hp{\le}\;+\frac{\de c}{p-1}\int_M\abs{DA}^2f^pH^{-2},
\end{aligned}
\end{equation}
\begin{equation}
\begin{aligned}
\al\int_MF^{-1}F^{ij}F_{;ij}f^p&=-\al \int_MF^{-1}F^{ij}F_if_jpf^{p-1}\\
&\msp[-100]-\al \int_MF^{-1}F^{ij,kl}h_{kl;j}F_if^p+\al \int_MF^{-2}F^{ij}F_iF_jf^p,
\end{aligned}
\end{equation}
which can be estimated by the right-hand side of \re{7.52}.
\begin{equation}
\begin{aligned}
-2\int_M\{h^{ij}-FnF^{ij}\}F^{-\al}F_{;ij}f^{p-1}&=\\
&\msp[-300] 2(p-1)\int_M\{h^{ij}-FnF^{ij}\}F_if_jf^{p-2}F^{-\al}\\
&\msp[-304]+2\int_M\{h^{ij}-FnF^{ij}\}_jF_iF^{-\al}f^{p-1}\\
&\msp[-304]-2\al\int_M\{h^{ij}-FnF^{ij}\}F_jF_iF^{-(1+\al)}f^{p-1}.
\end{aligned}
\end{equation}
In view of the estimate \re{7.49} the right-hand side of the preceding equality can be estimated as desired.

Finally, let us consider
\begin{equation}
\begin{aligned}
2(\al-1)\int_MF^{-1}F^{ij}F_if_jf^{p-1}&\le c\int_MF^{ij}f_if_jf^{p-2}\\
&\hp{\le}\;+c\int_M\abs {DA}^2H^{-2}f^p,
\end{aligned}
\end{equation}
which can be estimated as desired completing the proof of the lemma.
\ep

Now we can show that for large $p$ the $L^p$-norms of $f=f(t,\cdot)$ are uniformly bounded provided $\s$ is small enough.
\bl\lal{7.9}
Let $F$ be strictly concave, then there exist $C_1>0$ and $\s_0>0$ such that for all
\begin{equation}\lae{7.58}
p\ge c\e^{-2}\q\wed\q \s\le \min(\e^3p^{-\frac12}\frac1{4c},\s_0),
\end{equation}
where $c>1$ is the constant in \re{7.52}, the estimate
\begin{equation}
\norm{f}_{p.M}\le C_1\qq\A\, t\in[0,T^*)
\end{equation}
is valid, where $C_1=C_1(M_0)$ and $\s_0=\s_0(F,M_0)$.
\el
\bp
We multiply equation \re{7.34} with $pf^{p-1}$ and integrate by parts. Observing that the terms involving $K_N$ add up to be non-positive if $\s$ is small, $\s\le\s_0$, in view of \rl{7.4}, \re{7.17.1} and the fact that
\begin{equation}
1\le F^{kl}g_{kl}\le c_0,
\end{equation}
and by applying the estimate
\begin{equation}
F^{kl,rs}h_{kl;i}h_{rs;j}h^{ij}F^{-\al}\le -2\e^2\abs{DA}^2F^{-\al}
\end{equation}
which has already been used in the proof of \rl{7.5}, we obtain
\begin{equation}
\begin{aligned}
&\df {}t\int_Mf^p+p(p-1)\int_MF^{ij}f_if_jf^{p-2}+2\e^2p\int_M\abs{DA}^2F^{-\al}f^{p-1}\\
&\le \s p\int_MF^{ij}h_{ki}h^k_jf^p+\tfrac12 p(p-1)\int_MF^{ij}f_if_jf^{p-2}\\
&\hp{\le}\;+2\frac {cp}{p-1}\int_MF^{ij}F_iF_jf^{p-1}F^{-\al}\\
&\le \s p\int_MF^{ij}h_{ki}h^k_jf^p+\tfrac12 p(p-1)\int_MF^{ij}f_if_jf^{p-2}\\
&\hp{\le}\;+c\int_M\abs{DA}^2F^{-\al}f^{p-1},
\end{aligned}
\end{equation}
where we may choose $c$ to be the same constant that we used in \re{7.52}. Hence, we deduce, because of \re{7.58},
\begin{equation}\lae{7.63}
\begin{aligned}
\df {}t\int_Mf^p+\tfrac12p(p-1)\int_MF^{ij}f_if_jf^{p-2}+\e^2p\int_M\abs{DA}^2F^{-\al}f^{p-1}&\\
&\msp[-300]\le \s p\int_MF^{ij}h_{ki}h^k_jf^p.
\end{aligned}
\end{equation}
Choosing now 
\begin{equation}
\s\le \min(c_0^{-1}\e^3p^{-\frac12},\s_0),
\end{equation}
where $c_0>0$ will be specified below, and
\begin{equation}
\de=\e p^{-\frac12},
\end{equation}
we infer from \re{7.52} that the right-hand side of  inequality \re{7.63} can be estimated from above by
\begin{equation}
\begin{aligned}
\msp[40]\frac{\e p^{\frac12}}{c_0}\{\e^2\int_MF^{ij}h_{ki}h^k_jf^p\}&\le\\
&\msp[-210]\frac{\e p^{\frac12}}{c_0}\{\de^{-1}c(p-1)+c\}\int_MF^{ij}f_if_jf^{p-2}\\
&\msp[-214]+\frac{\e p^{\frac12}}{c_0}\{\de c(p-1)+c\}\int_M\abs{DA}^2F^{-\al}f^{p-1}\\
&\msp[-240]=c_0^{-1}\{p(p-1)c+\e p^\frac12 c\}\int_MF^{ij}f_if_jf^{p-2}\\
&\msp[-240]\hp{=}\;+c_0^{-1}\{\e^2(p-1)c+\e p^\frac12 c\}\int_M\abs{DA}^2F^{-\al}f^{p-1}\\
&\msp[-240]\le c_0^{-1}2 cp(p-1)\int_MF^{ij}f_if_jf^{p-2}\\
&\msp[-240]\hp{\le}\;+c_0^{-1}2c\e^2(p-1)\int_M\abs{DA}^2F^{-\al}f^{p-1},
\end{aligned}
\end{equation} 
where we again applied\re{7.58}.
\ep
Choosing
\begin{equation}
c_0=4c
\end{equation}
leads to
\begin{equation}
\df{}t\int_Mf^p\le 0\qq\A\, t\in [0,T^*)
\end{equation}
from which the result immediately follows.
\bt\lat{7.10}
Let $F$ be strictly concave or let $F=\tfrac1n H$ then there exist constants $\de>0$ and $c_0>0$ depending only on $F$ and $M_0$ such that
\begin{equation}\lae{7.69}
\abs{A}^2 -nF^2\le c_0F^{2-\de},
\end{equation}
or equivalently,
\begin{equation}\lae{7.70}
\abs A^2-\tfrac1n H^2\le c_0 H^{2-\de}.
\end{equation}
\et
\bp
When $F=\tfrac1n H$ we use the estimate \re{7.41.1} instead of \re{7.41} to obtain the result in \rl{7.9}. Then, in both cases, $F$  strictly concave or $F=\tfrac1n H$, the further arguments are essentially identical to those in Huisken's paper.
\ep
\br
In the proof of \rl{7.8}, \rl{7.9} and \rt{7.10} we used the fact that the sectional curvature $K_N$ satisfies
\begin{equation}
K_N\ge 0
\end{equation}
but only out of convenience. In case of the opposite sign slightly different arguments would have prevailed, since the terms stemming from the curvature of the ambient space are of lower order and can be handled fairly easily.
\er

Combining the estimate \re{7.70} with the regularity result of the rescaled hypersurfaces we shall prove that the rescaled hypersurfaces converge to a unit sphere in $C^\un(\Ss[n])$ exponentially fast provided $F$ is strictly concave or $F=\tfrac1n H$.
First, we prove:
\bl
Let $F$ be strictly concave or $F=\tfrac1n H$,  let $\tilde M(\tau)$ be the rescaled hypersurfaces and $\tilde h_{ij}$, $\tilde F$, etc.\ be the rescaled geometric quantities, then there are positive constants $c,$ $\de$ such that
\begin{equation}\lae{7.72}
\int_{\tilde M}\abs{D\tilde A}^2\le ce^{-\de\tau} \qq\A\, \tau_0\le \tau<\un,
\end{equation}
where
\begin{equation}
\tau_0=-\log\Theta(0,T^*),
\end{equation}
and where we emphasize that each geometric quantity is scaled separately by multiplying or dividing it  with appropriate powers of $\Theta$, and by pointing out that the scaled principal curvature are not the principal curvatures of $\tilde M$. This caveat applies especially to the integral in \re{7.72}.
\el
\bp
Consider the inequality \re{7.44}, where now $f$ is defined by choosing $\s=0$, i.e.,
\begin{equation}
f=F^{-2}(\abs A^2-nF^2).
\end{equation}
$f$ is scale invariant, hence we deduce from \re{7.69} and \frc{6.8} 
\begin{equation}\lae{7.75}
f\le c_0e^{-\de\tau}\qq\A\, \tau\ge \tau_0.
\end{equation}
All terms in inequality \re{7.44} scale like $f$, i.e., they are of order zero. Integrating over $M$, using integration by parts and rescaling the resulting inequality yields the result in view of \re{7.49} and \re{7.75}.
\ep
Applying now the interpolation inequalities for Sobolev norms, \cf \cite[Theorem 4.17]{adams}, we conclude  that $\abs{D\tilde A}$ decays exponentially fast in $C^\un(S^n)$, hence we conclude
\bl\lal{7.13}
There exists positive constants $c$, $\de$ such that
\begin{equation}\lae{7.76}
\tilde F_{\max}-\tilde F_{\min}\le c e^{-\de\tau}\qq\A\, \tau\ge\tau_0
\end{equation}
and
\begin{equation}\lae{8.68.1}
\norm{D\tilde F}\le c e^{-\de\tau}\qq\A\,\tau\ge \tau_0.
\end{equation}
\el
\bp
We first estimate the unscaled quantities in $M(t)$
\begin{equation}
\begin{aligned}
F_{\max}-F_{\min}\le \diam M(t)\sup_{M(t)}\norm{DF}\le c \diam M \sup_{M(t)}\abs{DA}
\end{aligned}
\end{equation}
to deduce
\begin{equation}
\tilde F_{\max}-\tilde F_{\min}\le c \diam \tilde M \sup_{\tilde M(\tau)}\abs{D\tilde A},
\end{equation}
hence the result. Note that
\begin{equation}
\abs{D\tilde A}^2=\Theta^2 g^{ij}h^k_{l;i}\Theta h^l_{k;j}\Theta
\end{equation}
and
\begin{equation}
\diam \tilde M=\diam M \Theta^{-1}\le c(\inf_{M(t)}\ka_1)^{-1} \Theta^{-1}\le\const,
\end{equation}
in view of Myers' theorem.
\ep

A similar lemma is also valid for the mean curvature:
\bl\lal{7.14}
There exists positive constants $c$, $\de$ such that
\begin{equation}\lae{7.81}
\tilde H_{\max}-\tilde H_{\min}\le c e^{-\de\tau}\qq\A\, \tau\ge\tau_0.
\end{equation}
\el
We are now ready to prove that the rescaled flow hypersurfaces converge to a sphere, to a geodesic sphere of radius $1$. 

\bl
Let $\abs{D\tilde u}$ be defined by
\begin{equation}
\abs{D\tilde u}^2=\s^{ij}\tilde u_i\tilde u_j,
\end{equation}
then there are positive constants $c$ and $\de$ such that
\begin{equation}\lae{8.98}
\abs{D\tilde u}\le ce^{-\de\tau}\qq\A\,\tau\ge \tau_0.
\end{equation}
\el
\bp
Let us look at the scaled scalar curvature equation in \fre{6.76}
\begin{equation}
\tilde u'=-\tilde F\Theta^{-1}\frac{\sin\Theta}{\cos\Theta} v+\tilde u.
\end{equation}
Define
\begin{equation}
\f=\log \tilde u
\end{equation}
and
\begin{equation}
w=\tfrac12 \abs{D\f}^2=\tfrac12\tilde u^{-2}\abs{D\tilde u}^2,
\end{equation}
then
\begin{equation}\lae{8.79}
\f'=- e^{-\f}\tilde F\Theta^{-1}\frac{\sin\Theta}{\cos\Theta}v+1,
\end{equation}
where we note that
\begin{equation}
v^2=1+\frac1{\sin^2u}\s^{ij}u_iu_j=1+\vt(u)^{-2}\s^{ij}\f_i\f_j
\end{equation}
\cf \re{7.91.1} and \fre{7.92.1}.

Differentiating now \re{8.79} with respect to $\f^kD_k$ we obtain
\begin{equation}
\begin{aligned}
w'=2 e^{-\f} w \tilde F\Theta^{-1}\frac{\sin\Theta}{\cos\Theta}v - e^{-\f}\tilde F\Theta^{-1}\frac{\sin\Theta}{\cos\Theta} v^{-1}\vt(u)^{-2}w_k\f^k+R,
\end{aligned}
\end{equation}
where $R$ decays exponentially in view of \re{8.68.1} or other more trivial estimates.  
The function
\begin{equation}
w_{\max} =\sup_{\tilde M(\tau)}w
\end{equation}
then satisfies
\begin{equation}
\begin{aligned}
w_{\max}'&=2 e^{-\f}w_{\max}\tilde F\Theta^{-1}\frac{\sin\Theta}{\cos\Theta} v+R\\
&\ge 2 e^{-\f}w_{\max}\tilde F\Theta^{-1}\frac{\sin\Theta}{\cos\Theta} v-c e^{-\de\tau}
\end{aligned}
\end{equation}
for almost every $\tau\ge\tau_0$, and we deduce
\begin{equation}\lae{8.107}
(w_{\max}-\tfrac  c{\de}e^{-\de\tau})'\ge 2 e^{-\f}w_{\max}\tilde F \Theta^{-1}\frac{\sin\Theta}{\cos\Theta}v.
\end{equation}
Hence,
\begin{equation}
\lim_{\tau\ra\un}w_{\max}
\end{equation}
exists and, because of
\begin{equation}
\un>\int_{\tau_0}^\un2 e^{-\f}w_{\max}\tilde F \Theta^{-1}\frac{\sin\Theta}{\cos\Theta}v\ge c\int_{\tau_0}^\un w_{\max},
\end{equation}
we obtain
\begin{equation}
\lim_{\tau\ra\un}w_{\max}=0,
\end{equation}
from which we conclude further, in view of \re{8.107},
\begin{equation}
w_{\max}(\tau)\le \tfrac c{\de}e^{-\de\tau}\qq\A\, \tau\ge \tau_0.
\end{equation}
\ep
As a corollary we can prove: 
\bc\lat{7.16}
The  rescaled flow hypersurfaces converge to the unit sphere in $C^\un(\Ss[n])$.
\ec
\bp
Let $\tilde u_k=\tilde u(\tau_k,\cdot)$ be a convergent subsequence in $C^\un(\Ss[n])$, then we deduce from \re{8.98} that the limit hypersurface is a sphere which  is the unit sphere, since the geodesic spheres with radius $\Theta$ intersect the hypersurfaces $M(t)=\graph  u$, \cf \fre{5.10}.  Since any convergent subsequence converges to the same limit, the corollary is proved.
\ep

Applying now the interpolation inequalities for the $C^m$-norms we can state:
\bt\lat{7.20}
Let $F$ be strictly concave or $F=\tfrac1nH$, then the rescaled function $\tilde u$ converges in $C^\un(\Ss[n])$ to the constant function $1$ exponentially fast. 
\et

Let us finally prove that the rescaled $F$-curvature converges to $1$ exponentially fast.
\bl
Let $F$ be strictly concave or $F=\tfrac1n H$, then
\begin{equation}
\lim_{t\ra T^*}F\Theta=1.
\end{equation}
\el
\bp
For fixed $0<t<T^*$ let
\begin{equation}
u(t,\xi_0)=u_{\max}(t).
\end{equation}
Then, by  applying the maximum principle, we infer that in that point
\begin{equation}
\ka_i\ge \frac{\cos u}{\sin u}
\end{equation}
and hence
\begin{equation}
\limsup_{t\ra T^*}\tilde F_{\max}\ge 1
\end{equation}
as well as
\begin{equation}
\liminf_{t\ra T^*}\tilde F_{\max}\ge 1.
\end{equation}
Looking at points $(t,\xi_0)$, where
\begin{equation}
u(t,\xi_0)=u_{\min}(t),
\end{equation}
we deduce the opposite inequalities for $\tilde F_{\min}$ proving the lemma, in view of the estimate \re{7.76}.
\ep
\bl\lal{7.18}
Let $F$ be strictly concave or $F=\tfrac1nH$, then there exist positive constants $c$, $\de$ such that
\begin{equation}
\abs{\tilde F(\tau,\cdot)-1}\le c e^{-\de\tau}\qq\A\, \tau\ge \tau_0.
\end{equation}
\el
\bp
We use the evolution equation for $F$. Let
\begin{equation}
\tilde F=F\Theta(t,T^*)
\end{equation}
and define
\begin{equation}
\tilde F_{\max}=F_{\max} \Theta(t,T^*),
\end{equation}
where
\begin{equation}
F_{\max}=\sup_{M(t)}F.
\end{equation}
Then we deduce from \fre{6.50}
\begin{equation}
\begin{aligned}
\tilde F_{\max}'&\le F^{ij}h_{ki}h^k_j F_{\max}\Theta^2\frac{\sin\Theta}{\cos\Theta}-\tilde F_{\max}\\
&\hp{\le}\;+K_NF^{ij}g_{ij}F_{\max}\Theta^2\frac{\sin\Theta}{\cos\Theta}
\end{aligned}
\end{equation}
for almost every $\tau\ge \tau_0$.

We now observe that
\begin{equation}
\abs{\Theta-\sin\Theta}\le c\Theta^2
\end{equation}
for small $\Theta$ and that
\begin{equation}
F^{ij}-\tfrac1n g^{ij}\le c(\sum_{i<j}(\ka_i-\ka_j)^2)^\frac12 g^{ij}
\end{equation}
\cf \re{7.49}, and, in view of \rl{7.4},
\begin{equation}
\tfrac1n\abs A^2-F^2\le c\sum_{i<j}(\ka_i-\ka_j)^2.
\end{equation}
Hence, we conclude
\begin{equation}
\tilde F_{\max}'\le (\tilde F_{\max}^2-1)\tilde F_{\max}+ce^{-\de\tau}\qq\A\,\tau\ge \tau_1,
\end{equation}
or equivalently,
\begin{equation}\lae{7.116}
(\tilde F_{\max}+\tfrac c\de e^{-\de\tau})'\le (\tilde F_{\max}^2-1)\tilde F_{\max}\qq\A\,\tau\ge\tau_1,
\end{equation}
where $\tau_1$ is sufficiently large such that $\Theta$ is small.

Suppose there exists $\tau_2>\tau_1$ such that
\begin{equation}\lae{7.117}
\tilde F_{\max}+\tfrac c\de e^{-\de\tau}<1
\end{equation}
in $\tau=\tau_2$, then this inequality is valid in a whole neighbourhood of $\tau_2$, since $\tilde F_{\max}$ is Lipschitz continuos, and we deduce from \re{7.116} that \re{7.117} is valid for all $\tau\ge \tau_2$ and
\begin{equation}
(\tilde F_{\max}+\tfrac c\de e^{-\de\tau})'\le0\qq\A\,\tau\ge \tau_2
\end{equation}
leading to the contradiction
\begin{equation}
1=\lim_{\tau\ra\un}(\tilde F_{\max}+\tfrac c\de e^{-\de\tau})\le \tilde F_{\max}(\tau_2)+\tfrac c\de e^{-\de\tau_2}<1.
\end{equation}
Thus, we conclude
\begin{equation}
\tilde F_{\max}-1\ge -\tfrac c\de e^{-\de\tau}\qq\A\,\tau\ge \tau_1.
\end{equation}
Defining
\begin{equation}
\tilde F_{\min}=\inf_{M(t)}F\Theta
\end{equation}
we deduce by an analogous argument
\begin{equation}
\tilde F_{\min}-1\le \tfrac c\de e^{-\de\tau}\qq\A\,\tau\ge \tau_1.
\end{equation}
Combining these two inequalities with inequality \re{7.76} completes the proof of the lemma.
\ep

\section{Inverse curvature flows}
Let the curvature functions $F$ govern the contracting curvature flows and their inverses $\tilde F$ the expanding flows
\begin{equation}
\dot x=\tilde F^{-1}\nu.
\end{equation}
A contracting flow converges to a point $x_0\in\Ss$ and are thus staying in the corresponding hemisphere $\mc H(x_0)$ for $t$ close to $T^*$, i.e., for $t_\de\le t<T^*$, and hence the corresponding expanding flow stays in the opposite hemisphere $\mc H(-x_0)$ for those values of $t$ and  converges to the equator. Since the flow is expanding, all flow hypersurfaces therefore stay in $\mc H(-x_0)$. The respective flow hypersurfaces are related by the Gau{\ss} map.

Fix a curvature $F$ to define a contracting flow and write the flow hypersurfaces $M(t)$ as graphs of a function $u$ with respect to geodesic polar coordinates centered in $x_0$ and write the polar hypersurfaces $M(t)^*$, which are the flow hypersurfaces of the corresponding inverse curvature flow, as graphs of a function $u^*$ with respect to geodesic polar coordinates centered in $-x_0$. This coordinate system will cover the  inverse curvature flow in the interval $[t_\de,T^*)$. Then we have:
\bl
The functions $u$, $u^*$ satisfy the relations
\begin{equation}\lae{9.2}
u_{\max}=\frac\pi2-u^*_{\min}\qq\A\, t\in [t_\de,T^*)
\end{equation}
and
\begin{equation}\lae{9.3}
u_{\min}=\frac\pi2-u^*_{\max}\qq\A\, t\in [t_\de,T^*).
\end{equation}
\el
\bp
Let $S_r(x_0)$ be a geodesic sphere around $x_0$ of radius $r$ and
\begin{equation}
S^*_r(x_0)=S_{r^*}(-x_0)
\end{equation}
be the polar sphere, then
\begin{equation}
\frac{\cos r}{\sin r}=\frac{\sin r^*}{\cos r^*},
\end{equation}
hence
\begin{equation}
r=\frac\pi2-r^*.
\end{equation}
Since the polar sets of convex bodies $\hat M_i$, $i=1,2$, satisfy
\begin{equation}\lae{9.7}
\hat M_1\su \hat M_2\q\im\q \hat M_2^*\su \hat M_1^*,
\end{equation}
\cf \cite[Corollary 9.2.10]{cg:cp}, we immediately deduce the relations \re{9.2} and \re{9.3} from \re{9.7}. 
\ep

\bc
There exists a positive constant $c$ such that
\begin{equation}\lae{9.8}
c^{-1}\le w=(\frac\pi2-u^*)\Theta^{-1}\le c\qq\A\, t\in [t_\de,T^*).
\end{equation}
\ec

\bl
Let
\begin{equation}
\abs{Dw}^2=\s^{ij}w_iw_j,
\end{equation}
then there exists a positive constant such that
\begin{equation}\lae{9.10}
\abs{Dw}^2\le c\qq\A\, t\in[t_\de,T^*).
\end{equation}
\el
\bp
Let $v$ be defined by
\begin{equation}
v^2=1+\frac1{\sin^2u^*}\s^{ij}u^*_iu^*_j,
\end{equation}
then
\begin{equation}\lae{9.12}
v\le e^{\bar\ka (u^*_{\max}-u^*_{\min})}\qq\A\, t \in [t_\de,T^*),
\end{equation}
where $\bar\ka$ is a positive upper bound for the principal curvature of the slices $\{x^0=\const\}$ that intersect The flow hypersurfaces, \cf \cite[inequality (2.7.83)]{cg:cp}, hence we conclude that for fixed $t$
\begin{equation}
\begin{aligned}
\frac1{\sin^2u^*}\s^{ij}u^*_iu^*_j&=v^2-1\le e^{2\bar\ka(u^*_{\max}-u^*_{\min})}-1\\
&\le c \sup_{M(t)}\frac{\cos u^*}{\sin u^*}(u^*_{\max}-u^*_{\min})\le c \Theta^{2},
\end{aligned}
\end{equation}
in view of $\re{9.8}$, hence the result.
\ep

The inverse curvature flow exists in the interval $[0,T^*)$ and is smooth. In order to prove this, we choose a point $y_0\in \hat M_0^*$ as the center of a geodesic polar coordinate system, then this system covers the whole flow, since the flow hypersurfaces are boundaries of strictly convex bodies. We have $C^0$ and $C^1$-estimates, \cf \re{9.12}, as well as $C^2$-estimates. Furthermore, $\tilde F$ is strictly positive on compact subintervals of $[0,T^*)$, hence the flow is smooth on compact subintervals.

For the rescaling process we may therefore restrict our attention to the interval $[t_\de,T^*)$, where we can  write the flow hypersurfaces as graphs in the coordinate system centered at $-x_0$. For $u^*$ we have the estimates \re{9.8} and \re{9.10}. Using then similar arguments as in the proofs of \frl{6.6} and \frt{6.12} we  conclude:
\bt
Let $u^*$ represent an inverse curvature curvature flow in $\Ss$ in the geodesic polar coordinate system specified above, where the curvature function and its inverse are both monotone and concave, then $u^*$ converges to the constant function $\frac\pi2$ in $C^\un(\Ss[n])$ such that for any $m\in\N$ the estimate
\begin{equation}
\abs{\tfrac\pi2 -u^*}_{m,\Ss[n]}\le c_m \Theta\qq\A\, t\in [0,T^*),
\end{equation}
is valid. The rescaled functions
\begin{equation}
w=(\tfrac\pi2-u^*)\Theta^{-1}
\end{equation}
are uniformly bounded in $C^\un(\Ss[n])$. When the curvature function $F$ of the corresponding contracting flow is strictly concave, or when $F=\tfrac1n H$, then  $w(\tau,\cdot)$ converges in $C^\un(\Ss[n])$ to the constant function $1$ exponentially fast.
\et

%\backmatter
%\includepdf[pages=-]{/Users/claus/Documents/Scanned-Documents/}
\bibliographystyle{hamsplain}
%\bibliography{mrabbrev,publications}
\providecommand{\bysame}{\leavevmode\hbox to3em{\hrulefill}\thinspace}
\providecommand{\href}[2]{#2}

%\listoffigures

%\cleardoublepage

%\thispagestyle{empty}
%\closegraphsfile
\end{document}